\documentclass[11pt]{article}

\RequirePackage[OT1]{fontenc}
\RequirePackage[aap,amsthm,amsmath,natbib]{imsart}
\RequirePackage[dvips,colorlinks]{hyperref}
\usepackage{mathrsfs, amssymb, dsfont}

\startlocaldefs \numberwithin{equation}{section}
\theoremstyle{plain}
\newtheorem{cond}{Condition}[section]
\newtheorem{cor}{Corollary}[section]

\newtheorem{lem}{Lemma}[section]
\newtheorem{prop}{Proposition}[section]
\newtheorem{rem}{Remark}[section]
\newtheorem{thm}{Theorem}[section]
\endlocaldefs

\begin{document}

\begin{frontmatter}

\title{Interacting Agent Feedback Finance Model}

\runtitle{Interacting Agent Feedback Finance Model}

\begin{aug}
\author{\fnms{Biao Wu} \ead[label=e1]{biaowu@math.mcmaster.ca}}

\runauthor{Biao Wu}

\affiliation{McMaster University}

\address{Department of Mathematics and Statistics\\
McMaster University\\
Hamilton Hall\\
1280 Main Street West\\
Hamilton, Ontario\\
Canada L8S 4K1 \printead{e1}}
\end{aug}

\begin{abstract}
We consider a financial market model which consists of a financial
asset and a large number of interacting agents classified into many
types. Different types of agents are heterogeneous in their price
expectations. Each agent can change its type based on the current
empirical distribution of the types and the equilibrium price, and
the equilibrium price follows a recursive price mechanism based on
the previous price and the current empirical distribution of the
types. The interaction among the agents, and the interaction between
the agents and the equilibrium price, feedback, are modeled. We
analyze the asymptotic behavior of the empirical distribution of the
types and the equilibrium price when the number of agents goes to
infinity. We give a case study of a simple example, and also
investigate the fixed points of empirical distribution and
equilibrium price of the example.
\end{abstract}



\end{frontmatter}

\section{Introduction}\label{se:Introduction}

Stochastic models of interacting systems play an important role in
population biology and statistical physics, c.f. \cite{Kingman} and
\cite{ChenMF}. In recent years a number of leading thinkers have
expressed the need for developing economic models that incorporate
interactions between agents and evolutionary mechanisms, see
 \cite{Palmer} and \cite{Farmer}. Agent-based models (ABMs), which arise from
many areas of science and social sciences such as ecology,
artificial intelligence, communication networks, sociology,
economics (see \cite{Ferber99}, \cite{Ferber95}, \cite{Ouelhadj96},
\cite{Wooldridge95}) are the ideal choice to attain this goal. The
following features of ABMs (see \cite{Goldstone}) are fundamental.
Firstly, precise mathematical formulation can be described by ABMs,
which make clear, quantitative and objective predictions possible.
Secondly, the explanations that link the analysis of the individual
agent level and the analysis of the emergent aggregate level can be
bridged by ABMs.

In this paper we will concentrate on the agent based modeling in a
financial market. At first, we give an account of some aspects of
the related works done by other authors. Black \cite{Black86}
classified traders as information traders and noise traders.
F\"{o}llmer and Schweizer \cite{Follmer&Schweizer} considered an
interacting agent financial model in which they used Black's
classification of traders. They assumed the number of traders being
countable and introduced an individual excess demand function which
takes a simple log-linear form. Lux (see \cite{Lux1995},
\cite{Lux1997} and \cite{Lux1998}) assumed the number of agents
being finite, and divided the traders into three types: chartists,
optimists, and pessimists. In Lux's model, the types of traders can
be changed probabilistically, according to the profitability of each
type; there are new entrants into the market and exits of current
traders from the market. Horst (see \cite{Horst2000},
\cite{Horst2001}, \cite{Horst2002} and \cite{Horst2005}) considered
the interacting agent models with local and global interactions.
Horst classified the traders into three types: fundamental
optimistic traders, fundamental pessimistic traders and noise
traders. Let $\mathcal{A}$ be a countable set of active agents and
$\mathcal{A}_n$ be a sequence of finite sets satisfying $\lim_{n\to
\infty}\mathcal{A}_n=\mathcal{A}$. At each period $t\ge 0$, each
fundamental trader has its mood, e.g., $x^a_t=+1$ for being
optimistic or $x^a_t=-1$ for being pessimistic. Let $C$ be a fixed
set of individual states, i.e., $x^a_t\in C$ for each $a\in
\mathcal{A}$ and $t\ge 0$. Let $x_t=\{x^a_t\}_{a\in \mathcal{A}}$.
Horst defined the empirical distribution as follows:
\begin{equation}\label{EmpricalDistHorst}
\rho_t=\rho(x_t)=\lim_{n\to\infty}\frac{1}{|\mathcal{A}_n|}\sum_{a\in
\mathcal{A}_n}\delta_{x^a_t}(\cdot).
\end{equation} $\rho(x_t)$ is called the mood of the market at time
$t$. The market mood drives the market price in the following way:
\begin{equation}\label{HorstRecu}
p_{t+1}=G(\rho_{t+1},p_t),
\end{equation} where $p_t$ is the market
price in period $t$ and $G$ is a certain function. Assuming a simple
log-linear structure for the excess demand function, Horst got the
recursive log price formula of the following form: for each period
$t\ge 0$,
\begin{equation}\label{HorstlogRecu}
\log p_{t+1}=f(\rho_{t+1})\log p_t+ g(\rho_{t+1}).
\end{equation}
The mood for each individual $a\in \mathcal{A}$ evolves as follows:
\begin{equation}\label{HorstIndmood}
\pi_a(x^a_{t+1}=s|x_t, e_{t+1},h_t)=\pi(x^a_{t+1}=s|x^a_t,
e_{t+1},h_t),
\end{equation} where $s\in C$, and $e_{t+1}\sim
Q(\rho(x_t);\cdot)$ is the signal of the market mood $\rho(x_t)$ and
$h_t$ is some (exogenous) economic fundamentals revealed in period
$t$. Therefore, in Horst's models, traders can change their types
during the evolution of the models and there are interactions among
traders. But there is no feedback of the price $p_t$ on the
evolution of individual state, i.e., the current market price has no
impact on the change of types for the next time period.

As was illustrated above, the empirical distribution of the types of
agents can link the behavior of individual agent level, the emergent
laws of aggregate level, and the equilibrium price of certain financial
asset. This paper is the second attempt of a systematic study of the
interacting agent financial system. Another working paper of the
author focuses on the multiagent models evolving in time-varying and
random environment, see \cite{Wu2a}. We construct the interacting
agent feedback finance model (IAFFM) by using agent based modeling.
The most general assumptions about the mechanism of IAFFM are as
follows:
\begin{enumerate}
  \item The time is in nonnegative integer units, denoted by $k\ge 0$.

  \item There are fixed $N\ge 2$ agents in the financial market at all times. There are no entries of new agents into the market or
exits of current agents from the market.

  \item There is a financial asset in certain market whose price $S_N(k)$ varies with time $k$.

  \item Each agent has one and only one internal state from an internal system with $r\ge 2$ states denoted by $1$, $\cdots$, $r$.
  The internal system does not change with time $k$, i.e, at each time period $k\ge 0$, there is no new state added to it and
  no existing state removed from it. The agents are classified into $r$ types according to their internal states, and
  $\mathbf{n}^N(k)=(n^N_1(k), \cdots, n^N_r(k))$ is the distribution of these agents among the $r$ types. It follows that $n^N_1(k)+ \cdots+
n^N_r(k)=N$ by the third assumption, and $\frac{\mathbf{n}^N(k)}{N}$ is
the empirical distribution of the types of agents at $k\ge 0$.

\item There exists a log price mechanism for the financial asset. Let $Z_+=\{0, 1, \cdots\}$,
$K_N=\{N^{-1}\mathbf{\alpha}: \mathbf{\alpha}\in (Z_+)^r, \sum_{i=1}^{r}\alpha_i=N\}$, and $g_N$ be defined on $[0,\infty)\times K_N\times R$.
At time $k\ge 1$, $\tilde{q}^N(\frac{k}{N})=\log S_N(k)$ is determined by $\tilde{q}^N(\frac{k-1}{N})$ and $\frac{\mathbf{n}^N(k)}{N}$ through the following recursive formula:
 \begin{equation}\label{eq:RecLogPriceTdqNgN}
\tilde{q}^N(\frac{k}{N})=\tilde{q}^N(\frac{k-1}{N})+\frac{1}{N}g_N(\frac{k}{N},
\frac{{\mathbf n}^N(k)}{N},\tilde{q}^N(\frac{k-1}{N})).
\end{equation}

\item Assume that $P_N(\cdot, \cdot, \cdot) =(p_{N, i, j}(\cdot, \cdot, \cdot))_{r\times r}$ is a deterministic {\em stochastic-matrix} valued function
 defined on $Z_+\times K_N\times R$ which represents the external environment of the multiagent system.

  \item Based on all the information of the agents' types, the equilibrium prices of the financial asset, and the external environment up to time $k$,
  each agent has an independent strategy of probabilistically choosing its type for the next time unit $k+1$.
  The strategy of an agent is realized by keeping or changing its type. The agents of the same type have a common strategy. That is to say,
  from time $k$ to $k+1$, the agents of type $i$ switch to type $j$ with probability $p_{N, i,
  j}(k, \frac{\mathbf{n}^N(k)}{N}, \tilde{q}^N(\frac{k}{N}))$. This process of changing types
  occurs locally among agents of the same type.

\end{enumerate}

Note that $S_N$ and $g_N$ depend on $N$, which means that the price mechanism does reflect the influence of
the size of the market. The gross performance of an economy consists of an external environment. The economic fundamentals of the financial asset can be reflected in $g_N$.
An example will be given in Section 4 to illustrate these. The assumptions 1-7 above will be used to mathematically formulate IAFFM.
 When we assume that $P_N(\cdot, \cdot, \cdot)$ or $g_N$ are random, we can construct IAFFM evolving in random environment.

  In this paper we will mainly study the asymptotic behaviors of $\{ (\frac{\mathbf{n}^N([Nt])}{N},$ $\tilde{q}^N(\frac{[Nt]}{N})), t\ge 0\}$
   as $N\to \infty$. One feature of our model is that the transition
structure of the types and log equilibrium price is time-inhomogeneous. Another feature is that the equilibrium price of the financial asset has feedback
 on the transition of the agent's types, instead of just being driven by the empirical distribution of the agents' types.
 Therefore, we modeled two kinds of interactions: the interaction among the agents, and the interaction between
the agents and the equilibrium price.

This paper is organized as follows. In Section \ref{se:FormResults}, we formulate IAFFM and state
the main Theorem (Theorem \ref{thm:ConvIAFFM}). In Section \ref{se:Example}, we give a complete case study of a simple IAFFM example. We make assumptions
for this example, verify its conditions required by Theorem \ref{thm:ConvIAFFM}, and discuss its fixed point problem. We also make connections
with the classical stock price formula for this example. In Section \ref{se:ProofofThm}, we give the proof of Theorem \ref{thm:ConvIAFFM}. In Appendix A,
B, and C, we give the proof of Lemmas \ref{qxqtdxBd}, \ref{rem:verifyCondiEofCor3.6Chp2}, and \ref{lem:LimitPtBeingCont}, respectively.

\section{Formulation of IAFFM and Main Theorem}\label{se:FormResults}

Let $R_+=[0, \infty)$ and $K=\{\mathbf{\alpha}: \mathbf{\alpha}\in (R_+)^r, \sum_{i=1}^{r}\alpha_i=1\}$.

{\cond\label{cond:regA} $A(t, \mathbf{x}, q)=(a_{i,j}(t, \mathbf{x}, q))_{r\times r}$ is a $r\times r$ matrix-valued
$c\grave{a}dl\grave{a}g$ function on $[0, \infty)\times K\times R$, which satisfies the conditions as follows: for each $(\mathbf{x},q)\in K\times R$,
\begin{itemize}
\item[1)] $A(\cdot, \mathbf{x}, q)\in D_{R^{r\times
r}}[0,\infty)$;
\item[2)] $A(t, \mathbf{x}, q) \mathbf{e}={\mathbf 0}_{r\times 1}$ for $t\ge 0$, where $\mathbf{e}=[1,\cdots,1]'$;
\item[3)] $a_{i,j}(t,\mathbf{x},q)\ge 0$ for $t\ge 0$, and $1\le
i,j\le r$, $i\not=j$.
\end{itemize}
}
{\cond\label{cond:regAN} For each $N\ge 1$, $A_N(t,
\mathbf{x}, q)=(a_{N,i,j}(t, \mathbf{x}, q))_{r\times r}$ is a
$r\times r$ matrix-valued function on $[0, \infty)\times K_N\times
R$, which satisfies the following slightly different conditions: for each $(\mathbf{x},q)\in K_N\times R$,
\begin{itemize}
\item[1)]  $A_N(\cdot, \mathbf{x}, q)\in D_{R^{r\times
r}}[0,\infty)$;
\item[2)] $A_N(t, \mathbf{x}, q) \mathbf{e}={\mathbf 0}_{r\times 1}$, $t\ge 0$;
\item[3)] $a_{N,i,j}(t,\mathbf{x},q)\ge 0$ for $t\ge 0$ and $1\le i,j\le r$, $i\not=j$.
\item[4)] $A_N(t, \mathbf{x}, q)$ is a constant on $t\in [\frac{k}{N}, \frac{k+1}{N})$ for $k\ge 0$.
\end{itemize}
}

Next, we specify the functions which are related to the log price
mechanism. For each $N\ge 1$, we have defined $g_N(t, \mathbf{x}, q)$ on $[0, \infty)\times
K_N\times R$. We also define a real valued function $g(t, \mathbf{x}, q)$ on $[0, \infty)\times K\times R$.
We are only concerned with functions which are `linearly' growing in $q$, i.e.
\begin{align}
\label{eq:FormOfg}&\,g(t,\mathbf{x}, q)=\varphi(t,\mathbf{x},
q)q+\psi(t,\mathbf{x}, q), \\
\label{eq:FormOfgN}&\, g_N(t,\mathbf{x}, q)=\varphi_N(t,\mathbf{x},
q)q+\psi_N(t,\mathbf{x}, q),
\end{align}
where $\varphi(t,\mathbf{x}, q)$ and $\psi(t,\mathbf{x}, q)$ are
functions on $[0, \infty)\times K\times R$, and
$\varphi_N(t,\mathbf{x}, q)$ and $\psi_N(t,\mathbf{x}, q)$ are
functions on $[0, \infty)\times K_N\times R$.

Let $\mathbf{X}^N(t)=\frac{{\mathbf n}^N([Nt])}{N}$,
$q^N(t)=\tilde{q}^N(\frac{[Nt]}{N})$ and
$\mathbf{Y}^N(t)=(\mathbf{X}^N(t),q^N(t))$.
(\ref{eq:RecLogPriceTdqNgN}) represented by the new notations is for
$t\ge \frac{1}{N}$,
\begin{equation}\label{eq:RecLogPriceqNtgN}
q^N(t)=q^N(t-\frac{1}{N})+\frac{1}{N}g_N(\frac{[Nt]}{N},
\mathbf{X}^N(t),q^N(t-\frac{1}{N})).
\end{equation}

Now we formulate the transition structure of IAFFM as follows. For $k\ge 0$, $(\mathbf{x},q)\in K_N\times R$, let
\begin{equation}\label{eq:IAFFMPNANk}
P_{N}[k, \mathbf{x}, q]=I+\frac{1}{N}A_N(\frac{k}{N}, \mathbf{x}, q),
\end{equation}
where $I$ is the identity matrix of order $r$. For each fixed $k\ge
0$, $(\mathbf{x},q)\in K_N\times R$, $P_{N}[k, \mathbf{x}, q]$ is a
stochastic matrix for large enough $N$. Assume that ${\mathbf
n}^N(0)$ and $\tilde{q}^N(0)$ are given. We define the time
inhomogeneous Markov chain $\{({\mathbf
n}^N(k),\tilde{q}^N(\frac{k}{N}))\}_{k=0}^{\infty}$ by mathematical
induction. Assume that for $k\ge 0$, $({\mathbf
n}^N(k),\tilde{q}^N(\frac{k}{N}))$ is defined or given, we want to
define $({\mathbf n}^N(k+1),\tilde{q}^N(\frac{k+1}{N}))$. By the
assumption 7 in Section \ref{se:Introduction}, for $1\le i\le r$,
each agent of type $i$ can change its type to $j$, with probability
$p_{N, i, j}(k, \mathbf{X}^N(\frac{k}{N}),\tilde{q}^N(\frac{k}{N}))$
($1\le j\le r$). Since the $n^N_{i}(k)$ agents of type $i$
independently make their transitions, the distribution of these
$n^N_{i}(k)$ agents among the $r$ types at time $k+1$ is a random
vector denoted by $\mathbf{\Xi}_{N, k, i}=(\xi_{N, k, i,1},\cdots,
\xi_{N, k, i,r})$ on some probability space, which satisfies
\begin{equation}\label{eq:IAFFMXiNkj}
\mathbf{\Xi}_{N, k, i}
\sim \mbox{multinomial}\biggl(n^N_i(k),
P_{N,i\cdot}[k,\mathbf{X}^N(\frac{k}{N}),\tilde{q}^N(\frac{k}{N})]\biggl),
\end{equation}
where $P_{N,i\cdot}[k,\mathbf{X}^N(\frac{k}{N}),\tilde{q}^N(\frac{k}{N})]$ is
the $i$-th row of matrix
$P_{N}[k,\mathbf{X}^N(\frac{k}{N}),\tilde{q}^N(\frac{k}{N})]$. Since the agents with different type change
their types independently, $\mathbf{\Xi}_{N,k,1}$, $\cdots$, $\mathbf{\Xi}_{N,k,r}$
are independent. Then, we define
\begin{equation}\label{IAFFMeq:nk+In}
{\mathbf
n}^N(k+1)\equiv\mathbf{\Xi}_{N,k,1}+\cdots+\mathbf{\Xi}_{N,k,r},
\end{equation}
and $\tilde{q}^N(\frac{k+1}{N})$ by formula
(\ref{eq:RecLogPriceTdqNgN}).

For each $N\ge 1$, we define the transition operators of
$\{(\mathbf{X}^N(\frac{k}{N}),\tilde{q}^N(\frac{k}{N})),k\ge 0\}$ as
follows: for each $k\ge 0$ and $f\in \bar{C}(K_N\times R)$, the set
of bounded continuous function on $K_N\times R$, for $(\mathbf{x},q)\in
K_N\times R$,
\begin{equation}\label{eq:TranOpYNIMMF}
S_{N, k}f(\mathbf{x},q)=E[f(\mathbf{X}^N(\frac{k+1}{N}),\tilde{q}^N(\frac{k+1}{N}
))|\mathbf{X}^N(\frac{k}{N})=\mathbf{x}, \tilde{q}^N(\frac{k}{N})=q ].
\end{equation}

We expect that the transpose $\mathbf{Y}'$ of any limit $\mathbf{Y}$ of $\mathbf{Y}^N$ is
a solution of the following differential equations:
\begin{align}
\label{eq:IAFFMlimitdiffFormx}\frac{d\mathbf{x}(t)}{dt}&=A(t,\mathbf{x}(t),q(t))' \mathbf{x}(t),  \\
\label{eq:IAFFMlimitdiffFormq}\frac{dq(t)}{dt}&=g(t,\mathbf{x}(t),q(t)),
\end{align}
which satisfy the initial conditions $\mathbf{x}(0)=\mathbf{X}(0)'$
and $q(0)=Q(0)$.

Let $K\times R$ be the state space for the limit process and denote by
$\bar{C}(K\times R)$ the set of bounded continuous functions on
$K\times R$, and define $C_c(K\times R)$, $C^1_c(K\times R)$,
$C^2_c(K\times R)$
by
\begin{equation}\label{eq:CKRIAFFMcompactSupp}
C_c(K\times R)=\{f\in \bar{C}(K\times R), f \mbox{ has compact
support on } K\times R\},
\end{equation}
 \begin{equation}\label{eq:C1KRIAFFM}
 C^1_c(K\times R)=\{f\in
C_c(K\times R), \frac{\partial f}{\partial {\mathbf x}} \mbox{ and }
\frac{\partial f}{\partial q}\mbox{ are continuous on } K\times R\},
\end{equation}
 \begin{equation}\label{eq:C2KRIAFFM}
 C^2_c(K\times R)=\{f\in C^1_c(K\times R),
\frac{\partial ^2 f}{\partial {\mathbf x}^2} \mbox{, }
\frac{\partial ^2 f}{\partial {\mathbf x}\partial q}\mbox{, and }
\frac{\partial ^2 f}{\partial q^2} \mbox{ are continuous on }
K\times R\}.
\end{equation}
We define the time-dependent generator
    $\{G_A(s), 0\le s<\infty\} $ on $C^1_c(K\times R)$: for each $f\in C^1_c(K\times R)$, and $s\ge 0$,
\begin{equation}\label{eq:GeneratorIAFFM}
G_A(s)f(\mathbf{x},q)=\mathbf{x}A(s,\mathbf{x},q )\frac{\partial
f}{\partial \mathbf{x}}'+g(s,\mathbf{x},q )\frac{\partial
f}{\partial q}, \hspace{5mm}(\mathbf{x},q)\in K\times R.
\end{equation}
$\mathscr{D}(G_A)= C^1_c(K\times R)$ is the common domain of
 the generator $\{G_A(s), 0\le s<\infty\} $, and $D= C^2_c(K\times R)$, is a subalgebra contained in $\mathscr{D}(G_A)$.

{\cond\label{cond:UniformConvForAN} For each $f\in C^2_c(K\times R)$ and $T>0$ there exist measurable
sets $\{F_N\}\subset R$ such that
\begin{align}
\label{eq:ANConvAUniOnF} \lim_{N\to \infty}\sup_{q\in F_N}&\sup_{\mathbf{x}\in
K_N}d_U(A_N(\cdot,\mathbf{x},q),A(\cdot,\mathbf{x},q))=0, \\
\label{eq:YnFnConInProbOneIFMM}
\lim_{N\to \infty}&P\{q^N(t)\in F_N, 0\le t\le T\}=1,
\end{align}
where $d_U$ is the uniform metric on $D_{R^{r\times r}}[0, \infty)$ defined by
\begin{equation*}
d_U(\mathbf{u},\mathbf{v})=\int_0^{\infty} e^{-s}\sup_{0\le t\le
s}[\|\mathbf{u}(t)-\mathbf{v}(t)\|_r \wedge 1]ds, \mbox{
$\mathbf{u}$, $\mathbf{v}\in D_{R^{r\times r}}[0, \infty)$},
\end{equation*}
and  $\|\cdot\|$ is the matrix norm.
}

{\cond\label{cond:UniformConvForVarphiPsi}  For any $T>0$,
\begin{align}
\label{VarphiNconvToVarphi}&\,\lim_{N\to \infty}\sup_{0\le t\le
T}\sup_{(\mathbf{x},q)\in
 K_N\times R}|\varphi_N(t,\mathbf{x},q)-\varphi(t,\mathbf{x},q)|=0,
 \\
\label{PsiNconvToPsi}&\,\lim_{N\to \infty}\sup_{0\le t\le
T}\sup_{(\mathbf{x},q)\in
 K_N\times R}|\psi_N(t,\mathbf{x},q)-\psi(t,\mathbf{x},q)|=0.
\end{align}
}

{\cond\label{cond:ExpGrowth} For any compact subset $\tilde{K}$
of $K\times R$, there exist $C>0$ and $\lambda>0$, such that $A(t,
\mathbf{x}, q)$, $\varphi(t, \mathbf{x}, q)$ and $\psi(t,
\mathbf{x}, q)$ satisfy that
 \begin{align}
\label{eq:IAFFMA(s)ExpGrowth} \sup_{(\mathbf{x},q)\in \tilde{K}}\|A(t,\mathbf{x},q)\| \le &\, C e^{\lambda t},   \\
\label{eq:IAFFMVarphi(s)ExpGrowth} \sup_{(\mathbf{x},q)\in
\tilde{K}}|\varphi(t,\mathbf{x},q)| \le &\, C e^{\lambda t},  \\
\label{eq:IAFFMPsi(s)ExpGrowth} \sup_{(\mathbf{x},q)\in
\tilde{K}}|\psi(t,\mathbf{x},q)| \le &\, C e^{\lambda t}.
\end{align}
}

We introduce the following notations:
\begin{align}
\label{defnbiForOneTor} b_i(t,\mathbf{y})=&\,\sum_{j=1}^r x_j a_{j
i}(t,\mathbf{y}), \mbox{ $1\le i\le r$, and $\mathbf{y}=(\mathbf{x}, q)$ } \\
\label{defnbr+1} b_{r+1}(t,\mathbf{y})=&\,g(t, \mathbf{y}),
\end{align}
\begin{equation}\label{eq:BigBtHatBt}
\mathbf{b}(t, \mathbf{y})=(b_1(t,\mathbf{y}),\cdots,
b_{r+1}(t,\mathbf{y}))\mbox{ and }\hat{\mathbf{b}}(t,
\mathbf{y})=(b_1(t,\mathbf{y}),\cdots, b_{r}(t,\mathbf{y})).
\end{equation}
 It is clear that (\ref{eq:IAFFMlimitdiffFormx}) and
(\ref{eq:IAFFMlimitdiffFormq}) with the initial conditions
$\mathbf{x}(0)=\mathbf{X}(0)$ and $q(0)=Q(0)$ are equivalent to the
integral equations
\begin{equation}\label{eq:DetermInteEq}
{\mathbf{y}(t)=\mathbf{y}(0)+\int_0^t \mathbf{b}(s,\mathbf{y}(s))'
ds, \mbox{ }t\ge 0}
\end{equation}
or
\begin{align}
\label{eq:IAFFMlimitIntegFormx}\mathbf{x}(t)&=\mathbf{x}(0)'+\int_0^t \hat{\mathbf{b}}(s,\mathbf{x}(s),q(s))'ds,  \\
\label{eq:IAFFMlimitIntegFormq}q(t)&=q(0)+\int_0^t
g(s,\mathbf{x}(s),q(s))ds.
\end{align} where $\mathbf{y}(0)=(\mathbf{x}(0),q(0))$.
(\ref{eq:IAFFMlimitIntegFormx}) and (\ref{eq:IAFFMlimitIntegFormq}) are {\em nonlinear Volterra
equations of the second kind}.

(\ref{eq:IAFFMlimitIntegFormx}) and (\ref{eq:IAFFMlimitIntegFormq}) are assumed to satisfy the semi-Lipschitz conditions, which are guaranteed by the
following condition. The definition and property of semi-Lipschitz condition are included in Section \ref{se:ProofofThm}.
{\cond\label{cond:semiLipForIAFFM} For any $T>0$, there exists a nonnegative $\mathbb{L}^2([0,T], R)$ function
$C_T(t)$, such that for any $0\le t\le T$ and $\mathbf{y}_1=(\mathbf{x}_1,q_1), \mathbf{y}_2=(\mathbf{x}_2, q_2) \in K\times R$,
\begin{equation}\label{LipschitzForHatBt}
\|\hat{\mathbf{b}}(t,\mathbf{y}_1)-\hat{\mathbf{b}}(t,\mathbf{y}_2)\|\le
C_T(t)\|\mathbf{y}_1-\mathbf{y}_2\|.
\end{equation}
$g$ satisfies that for fixed $T>0$ and any $K$-valued continuous
function $\mathbf{x}(t)$ on $[0,T]$, (\ref{eq:IAFFMlimitIntegFormq})
determines a unique solution $q_{\mathbf{x}}(t)$ on $[0,T]$, and for
any two solutions $(\mathbf{x}(t), q_{\mathbf{x}}(t))$,
$(\tilde{\mathbf{x}}(t),q_{\tilde{\mathbf{x}}}(t))$ of
 (\ref{eq:IAFFMlimitIntegFormx}) and (\ref{eq:IAFFMlimitIntegFormq})
 on $[0,T]$, there exists a nonnegative
measurable function $B(s)$ on $[0,T]$, such that $q_{\mathbf{x}}(t)$
and $q_{\tilde{\mathbf{x}}}(t)$ satisfy
\begin{equation}\label{DiffOfTwoSolnsforq}
\int_0^s
 (q_{\mathbf{x}}(t)-q_{\tilde{\mathbf{x}}}(t))^2 dt\le B(s)
 \int_0^s\|\mathbf{x}(t)-\tilde{\mathbf{x}}(t)\|^2 dt \mbox{ for }
 0\le s\le T
 \end{equation}
 and
\begin{equation}\label{CTtBtInteg}
\int_0^T [\int_0^s C_T(t)^2 dt (1+B(s))]ds
 <\infty.
 \end{equation}
}

{\thm\label{thm:ConvIAFFM} Assume that $A(t, \mathbf{x},
q)=(a_{i,j}(t, \mathbf{x}, q))_{r\times r}$ satisfies the condition
\ref{cond:regA} and $\{A_N(t,\mathbf{x},q)=(a_{N,i,j}(t, \mathbf{x},
q))_{r\times r}\}$ satisfies the condition \ref{cond:regAN}. Assume
that $g(t,\mathbf{x},q)$ and $\{g_N(t,\mathbf{x},q)\}$ satisfy
(\ref{eq:FormOfg}) and (\ref{eq:FormOfgN}). Assume that $\varphi$,
$\psi$, $\frac{\partial \varphi}{\partial {\mathbf x}}$,
$\frac{\partial \psi}{\partial {\mathbf x}}$ are continuous. Assume
that the condition \ref{cond:UniformConvForAN} holds for $\{A_N\}$
$\{q^N\}$, and $A$, and the condition
\ref{cond:UniformConvForVarphiPsi} holds for $\{ \varphi_N\}$,
$\varphi$ and $\{\psi_N\}$, $\psi$. Assume that $A(t, \mathbf{x},
q)$, $\varphi(t, \mathbf{x}, q)$ and $\psi(t, \mathbf{x}, q)$
satisfy the condition \ref{cond:ExpGrowth} and are bounded on
$[0,T]\times K\times R$ for any $T>0$. Define $
\mathbf{b}(t,\mathbf{y})$ and $ \hat{\mathbf{b}}(t,\mathbf{y})$ by
(\ref{defnbiForOneTor}), (\ref{defnbr+1}) and (\ref{eq:BigBtHatBt}).
Assume that either $\mathbf{b}(t,\mathbf{y})$ satisfies the
Lipschitz condition or the condition \ref{cond:semiLipForIAFFM}
holds for (\ref{eq:IAFFMlimitIntegFormx}) and
(\ref{eq:IAFFMlimitIntegFormq}). Suppose that
$P(\mathbf{Y}^N(0))^{-1}\Rightarrow\mu$ for some
$\mu\in\mathscr{P}(K\times R)$, then there exists on some
probability space $(\Omega, \mathscr{F}, P)$ a stochastic process
$\mathbf{Y}$ satisfying $P(\mathbf{Y}(0))^{-1}=\mu$, which is the
unique solution of the $C_{K\times R}[0,\infty)$-martingale problem
for $(G_A, \mu)$ restricted to $C^2_c(K\times R)$ and for $P$-a.s.
$\omega\in \Omega$, $\mathbf{Y}(t,\omega)'$ is the unique solution
of the integral equations (\ref{eq:IAFFMlimitIntegFormx}) and
(\ref{eq:IAFFMlimitIntegFormq}), such that $\mathbf{Y}^N\Rightarrow
\mathbf{Y}$. }

\newpage

\section{Case Study of a Simple Example}\label{se:Example}

\subsection{A simple example of IAFFM}\label{sse:AGMEOIAFFM}

There are $N$ agents in a financial market which consist of
fundamentalists, optimists, and pessimists, see \cite{Lux1998}. At each
time $k\ge 0$, the number of fundamentalists, optimists and
pessimists are $n^N_1(k)$, $n^N_2(k)$, and $n^N_3(k)$ respectively.
This interacting agent system is closed, i.e., there are no new
entrants into the market or quits of current traders from the
market. Therefore $n^N_1(k)+n^N_2(k)+n^N_3(k)=N$ for all $k\ge 0$,
and $\mathbf{n}^N(k)=(n^N_1(k), n^N_2(k), n^N_3(k))$ is the vector
of types.

Each type of agents have their own excess demand functions. The
excess demand functions have the form of F\"{o}llmer and Schweizer,
see \cite{Follmer&Schweizer}. Assume that there exists a probability
space $(\Omega, \mathscr{F}, P)$. At time unit $k \ge 1$, for each $
\omega \in \Omega$, and a proposed price $p$, each agent $a$ has an
excess demand function $e_{a}^N(k,p,\omega)$ which is given by
\begin{equation}\label{eq:classicfunc1}
e_{a}^N(k,p,\omega)=\alpha_{a}^N(\frac{k}{N},\omega)\log({\hat{S}_{a,k}^N(\omega)}/p)+\frac{\delta^N_{a}(\frac{k}{N},\omega)}{N}.
\end{equation}
Here $\delta^N_{a}(\frac{k}{N},\omega)$ can be viewed as the total
liquidity demand observed by agent $a$ at time $k\ge 0$ and
$\frac{\delta^N_{a}(\frac{k}{N},\omega)}{N}$ is the average
liquidity demand, and $\hat{S}_{a,k}^N(\omega)$ denotes an
individual reference level of agent $a$ at time $k$. The time scale
for $\alpha_{a}^N$ and $\delta_a^N$ is of $\frac{1}{N}$, instead of
1.

Denote the individual reference level $\hat{S}_{a,k}^N(\omega)$, the
coefficients $\alpha_{a}^N(\frac{k}{N},\omega)$ and
$\delta^N_{a}(\frac{k}{N},\omega)$  for fundamentalist, optimist,
and pessimist by $\hat{S}_{i}^N(k,\omega)$,
  $\alpha_{i}^N(\frac{k}{N},\omega)$ and $\delta^N_{i}(\frac{k}{N},\omega)$ ($1\le i\le 3$) respectively. We assume individual reference levels
  as follows:
\begin{align}
\label{eq:fundref} \log \hat{S}_{1}^N(k,\omega)=&\log S_N(k-1)
+\frac{\beta_{1}^N(\frac{k}{N},\omega)}{N}(\log S_N(k-1)- \log
F^N(\frac{k}{N},\omega)),  \\
\label{eq:optref} \log\hat{S}_{2}^N(k,\omega)=&\log S_N(k-1)
+\beta_{2}^N(\frac{k}{N},\omega)(\log S_N(k-1)- \log p) , \\
\label{eq:pessref}\log \hat{S}_{3}^N(k,\omega)=&\log S_N(k-1)
+\beta_{3}^N(\frac{k}{N},\omega)(\log S_N(k-1)- \log p),
\end{align}
where random coefficients $\beta_{i}^N(\frac{k}{N},\omega) \le 0$
($1\le i\le 3$), $F^N(\frac{k}{N},\omega)$ is the fundamental value
of the asset at time $k$. Note that only
$\beta_{1}^N(\frac{k}{N},\omega)$ in the fundamentalists' reference
level is divided by $N$, the market size of the agents. One
reasonable explanation for it is that fundamentalists know the
evolution of fundamental value $F^N(\frac{k}{N},\omega)$ of the
asset, they value less in their excess demand functions the
difference between $\log S_N(k-1)$ and $\log
F^N(\frac{k}{N},\omega)$. Note also that the time scale for
$\beta_1^N$, $\beta_2^N$, $\beta_3^N$ and $\log F^N$ is also of
$\frac{1}{N}$.

For each $k\ge 1$, if we assume the price $S_N(k-1)$ and
$\mathbf{n}^N(k)$ are known, then the equilibrium log price $\log
S_N(k)$ is determined by the market clearing condition of zero
excess demand:
\begin{equation}
\sum_{a}e_{a}^N(k,S_N(k),\omega)=0,
\end{equation}
i.e.,
\begin{equation}\label{eq:equilibrium}
\sum_{i=1}^3 n^N_i(k)\biggl[\alpha_{i}^N(\frac{k}{N})(\log\hat{S}_{i}^N(k)-\log
S_N(k))+\frac{\delta^N_{i}(\frac{k}{N})}{N}\biggl]=0.
\end{equation}
We have omitted $\omega$ in the random variables in the above
equation and in the rest of this subsection. Substituting
(\ref{eq:fundref}), (\ref{eq:optref}), (\ref{eq:pessref}) into
(\ref{eq:equilibrium}), and solving for $\log p$ as $\log S_N(k)$,
we get
\begin{equation}\label{eq:reculogprice}
\begin{aligned}&\,\log S_N(k) -\log S_N(k-1) \\
= &\, \frac{1}{N}
\cdot\frac{n^N_1(k)\alpha_{1}^N(\frac{k}{N})\beta_{1}^N(\frac{k}{N})(\log
S_N(k-1)-\log F^N(\frac{k}{N}))}
{n^N_1(k)\alpha_{1}^N(\frac{k}{N})+n^N_2(k)\alpha_{2}^N(\frac{k}{N})(1+\beta_{2}^N(\frac{k}{N}))+n^N_3(k)\alpha_{3}^N(\frac{k}{N})(1+\beta_{3}^N(\frac{k}{N}))} \\
                    &  +\frac{1}{N} \cdot \frac{n^N_1(k)\delta^N_{1}(\frac{k}{N})+n^N_2(k)\delta^N_{2}(\frac{k}{N})+n^N_3(k)\delta^N_{3}(\frac{k}{N})}
                    {n^N_1(k)\alpha_{1}^N(\frac{k}{N})+n^N_2(k)\alpha_{2}^N(\frac{k}{N})(1+\beta_{2}^N(\frac{k}{N}))+n^N_3(k)\alpha_{3}^N(\frac{k}{N})(1+\beta_{3}^N(\frac{k}{N}))}.
\end{aligned}
\end{equation}
(\ref{eq:reculogprice}) is a recursive log price formula.

Let $K^3_N =\{N^{-1}\mathbf{x}: \mathbf{x}\in (Z_+)^3, \sum_{i=1}^3
x_i=N\}$. Define $g_N$ on $[0,\infty)\times K^3_N\times R$ as
follows:
\begin{equation}\label{eq:gNfuncIAFMMexamp}
\begin{aligned}
g_N(t,\mathbf{x},q)
=&\,\frac{x_1
\alpha_1^N(t)\beta_1^N(t)q}{x_1 \alpha_1^N(t)+x_2
\alpha_2^N(t)(1+\beta_2^N(t))+x_3 \alpha_3^N(t)(1+\beta_3^N(t))}
\\ &\,+\frac{x_1 \delta^{N}_{1}(t)+x_2 \delta^{N}_{2}(t)+x_3
\delta^{N}_{3}(t)-x_1 \alpha_1^N(t)\beta_1^N(t)\log F^N(t)}{x_1
\alpha_1^N(t)+x_2 \alpha_2^N(t)(1+\beta_2^N(t))+x_3
\alpha_3^N(t)(1+\beta_3^N(t))}.
\end{aligned}
\end{equation}
Then (\ref{eq:reculogprice}) can be represented by
$g_N(t,\mathbf{x},q)$ as
\begin{equation}\label{eq:logpricegN}
\log S_N(k) = \log S_N(k-1)+\frac{1}{N}g_N(\frac{k}{N},
 \frac{\mathbf{n}^N(k)}{N}, \log S_N(k-1)).
\end{equation}
Note that $g_N(t,\mathbf{x},q)$ defined by
(\ref{eq:gNfuncIAFMMexamp}) is a random function. We have actually
made preparations for the IAFFM evolving in a random environment.

We also define $A_N(t, \mathbf{x}, q)=(a_{N,i,j}(t, \mathbf{x},
q))_{3\times 3}$ on $[0, \infty)\times K_N^3\times R$ satisfying the
Condition \ref{cond:regAN}. Then we
can specify the transition structure for $\{(\mathbf{n}^N(k),\log
S_N(k))\}_{k=0}^{\infty}$ same way as that in Section
\ref{se:FormResults}.

In this example, $\alpha_i^N$, $\beta_i^N$, $\delta_i^N$ ($1\le i\le
3$), and $F^N$ consist of the external environment of the
interacting agent financial system. $F^N$ is the economic
fundamental of the financial asset.

\subsection{Assumptions and verifications of the example}\label{sse:CAVOTGME}

In this subsection, we make assumptions for the example and verify
the conditions of Theorem \ref{thm:ConvIAFFM} for this example.

Assume that $\alpha_i^N(t)$, $\beta_i^N(t)$, and $\delta_i^N(t)$ for
$1\le i\le 3$ and $\log F^N(t)$ are real valued functions defined on
$[0,\infty)$. Assume also that for
any $\mathbf{x}\in K_N^3$ and $t\ge 0$, $x_1 \alpha_1^N(t)+x_2
\alpha_2^N(t)(1+\beta_2^N(t))+x_3
\alpha_3^N(t)(1+\beta_3^N(t))\not=0$, which justify the definition of $g_N$ in (\ref{eq:gNfuncIAFMMexamp}). 
Then we can define real valued functions $\varphi_N$, $\psi_N$ on
$[0,\infty)\times K^3_N $ as follows:
\begin{equation}\label{varphiNExmp}
\varphi_N(t,\mathbf{x})=\frac{x_1 \alpha_1^N(t)\beta_1^N(t)}{x_1
\alpha_1^N(t)+x_2 \alpha_2^N(t)(1+\beta_2^N(t))+x_3
\alpha_3^N(t)(1+\beta_3^N(t))},
\end{equation}
\begin{equation}\label{psiNExamp}
\psi_N(t,\mathbf{x})=\frac{x_1 \delta^N_{1}(t)+x_2
\delta^N_{2}(t)+x_3 \delta^N_{3}(t)-x_1
\alpha_1^N(t)\beta_1^N(t)\log F^N(t)}{x_1 \alpha_1^N(t)+x_2
\alpha_2^N(t)(1+\beta_2^N(t))+x_3 \alpha_3^N(t)(1+\beta_3^N(t))}.
\end{equation}
$g_N$, $\varphi_N$ and $\psi_N$ so defined satisfy the relation
(\ref{eq:FormOfgN}).

Now we make the general assumptions which justify the description of the limit process
$\{\mathbf{y}(t)=(\mathbf{x}(t),q(t))$, $0\le t<\infty\}$.

Let $K^3=\{\mathbf{x}: \mathbf{x}\in (R_+)^3, \sum_{i=1}^3
x_i=1\}$. Let $\alpha_i(t)$, $\beta_i(t)$, $\delta_i(t)$ ($1\le i\le 3$), and $ F(t)$ are continuous real valued functions
defined on $[0,\infty)$, where $F(t)>0$ for each $t\ge 0$. For any
$\mathbf{x}\in K^3$ and $t\ge 0$, $x_1 \alpha_1(t)+x_2
\alpha_2(t)(1+\beta_2(t))+x_3 \alpha_3(t)(1+\beta_3(t))\not= 0$.
$\alpha_i(t)$ and $\beta_i(t)$ ($1\le i\le 3$) are bounded on $[0, \infty)$. There exist constants $C>0$ and $\lambda>0$ such that for
$t\in [0, \infty)$,
\begin{equation}\label{BdExpGrowthDeltaF}
|\delta_i(t)|\le C e^{\lambda t} \mbox{ for $1\le i\le 3$ and } |\log F(t)|\le C e^{\lambda t}.
\end{equation}
$A(t, \mathbf{x},q)$ satisfies the Condition \ref{cond:regA} and is bounded on $[0, T]\times K^3
\times R$ for any $T>0$. For any compact subset $\hat{K}$ of $R$,
there exist $C>0$ and $\lambda>0$ such that for any $t\ge 0$,
\begin{equation}\label{BdIAFFMA(s)ExpGrowthexmp}
\sup_{(\mathbf{x},q)\in K^3\times \hat{K}}\|A(t,\mathbf{x},q)\| \le  C e^{\lambda t}.
\end{equation}
$A(t,\mathbf{y})$ also satisfies the Lipschitz condition.

Based on the above assumptions, we can define real valued functions $g$, $\varphi$, $\psi$ on
$[0,\infty)\times K^3\times R$ as follows:
\begin{equation}
\begin{aligned}\label{eq:gfuncIAFMMexampSimpleOne}
g(t,\mathbf{x},q)=&\,\frac{x_1
\alpha_1(t)\beta_1(t)q}{x_1 \alpha_1(t)+x_2
\alpha_2(t)(1+\beta_2(t))+x_3 \alpha_3(t)(1+\beta_3(t))}
\\ &\,+\frac{x_1 \delta_{1}(t)+x_2 \delta_{2}(t)+x_3
\delta_{3}(t)-x_1 \alpha_1(t)\beta_1(t)\log F(t)}{x_1
\alpha_1(t)+x_2 \alpha_2(t)(1+\beta_2(t))+x_3
\alpha_3(t)(1+\beta_3(t))}.
\end{aligned}
\end{equation}
\begin{equation}
\begin{aligned}\label{varphiExmp}
\varphi(t,\mathbf{x})=\frac{x_1 \alpha_1(t)\beta_1(t)}{x_1
\alpha_1(t)+x_2 \alpha_2(t)(1+\beta_2(t))+x_3
\alpha_3(t)(1+\beta_3(t))},
\end{aligned}
\end{equation}
\begin{equation}
\begin{aligned}\label{psiExamp}
\psi(t,\mathbf{x})=\frac{x_1 \delta_{1}(t)+x_2 \delta_{2}(t)+x_3
\delta_{3}(t)-x_1 \alpha_1(t)\beta_1(t)\log F(t)}{x_1
\alpha_1(t)+x_2 \alpha_2(t)(1+\beta_2(t))+x_3
\alpha_3(t)(1+\beta_3(t))}.
\end{aligned}
\end{equation}
Then $g$, $\varphi$ and $\psi$ satisfy
the relation (\ref{eq:FormOfg}).

{\rem\label{rem:CommOnGAssump} We make the following immediate comments based on the above general assumptions.
 \begin{description}
 \item[(1)] (\ref{BdIAFFMA(s)ExpGrowthexmp}) implies that (\ref{eq:IAFFMA(s)ExpGrowth}) holds for $A(t, \mathbf{x}, q)$.
\ref{BdExpGrowthDeltaF} implies that (\ref{eq:IAFFMVarphi(s)ExpGrowth}) holds for
$\varphi$ and (\ref{eq:IAFFMPsi(s)ExpGrowth}) holds for $\psi$. Therefore, the Condition \ref{cond:ExpGrowth} holds.
As a result, $\varphi(t, \mathbf{x})$ and $\psi(t,
\mathbf{x})$ are bounded for $(t,\mathbf{x}, q)\in [0, T]\times K^3
\times R$ for any $T>0$, viewing $q$ as a dummy variable of
$\varphi(t, \mathbf{x})$ and $\psi(t, \mathbf{x})$.

\item[(2)]The boundedness of $\alpha_i(t)$ ($1\le i\le 3$) on $[0, \infty)$ is just assumed for convenience, since the numerators and denominators of $g$,
$\varphi$, and $\psi$ are linear in $\alpha_i$ ($1\le i\le r$).

\item[(3)] (\ref{eq:IAFFMlimitIntegFormq}) is equivalent to
\begin{equation}\label{dqtdtgtxq}
\frac{dq(t)}{dt}=g(t,\mathbf{x}(t),q(t))
\end{equation}
with initial condition $q(0)=y_4(0)$. By the form of $g$, fix $T>0$, we know that for any
given $K^3$ valued function $\mathbf{x}(t)$ on $[0,T]$ , there
exists a unique solution $q_{\mathbf{x}}(t)$ of (\ref{dqtdtgtxq}) on
$[0, T]$.

\item[(4)] By Remark \ref{rem:LipschitzBigb},
$\hat{\mathbf{b}}(t,\mathbf{y})=\mathbf{x}A(t,\mathbf{y})$ satisfies
the Lipschitz condition since $A(t,\mathbf{y})$ is bounded on $[0,T]\times K^3\times R$ for any
$T>0$.

\end{description}
}

We prove in the next lemma that (\ref{DiffOfTwoSolnsforq}) and (\ref{CTtBtInteg}) are satisfied and
it follows the uniqueness of the $C_{K^3\times
R}[0,\infty)$-martingale problem for $(G_A, \mu)$.

{\lem\label{qxqtdxBd} Make the general assumptions on $\alpha_i(t)$,
$\beta_i(t)$, and $\delta_i(t)$ ($1\le i\le 3$), $ F(t)$ and $A(t,
\mathbf{x}, q)$ above. Then for fixed $T>0$ and any two solutions
$(\mathbf{x}(t), q_{\mathbf{x}}(t))$,
$(\tilde{\mathbf{x}}(t),q_{\tilde{\mathbf{x}}}(t))$ of
 (\ref{eq:IAFFMlimitIntegFormx}) and (\ref{eq:IAFFMlimitIntegFormq})
 on $[0,T]$, there exists $M>0$, such that
$q_{\mathbf{x}}(t)$ and $q_{\tilde{\mathbf{x}}}(t)$ satisfy
\begin{equation}\label{DiffOfTwoSolnsforqExamp}
\int_0^s (q_{\mathbf{x}}(t)-q_{\tilde{\mathbf{x}}}(t))^2 dt\le M
 \int_0^s\|\mathbf{x}(t)-\tilde{\mathbf{x}}(t)\|^2 dt \mbox{ for }
 0\le s\le T.
 \end{equation}
We conclude that the $C_{K^3\times R}[0,\infty)$-martingale problem for $(G_A, \mu)$
restricted to $C^2_c(K^3\times R)$ has at most one solution.
 }

We give in the next corollary the conditions which can
guarantee the weak convergence of $\{\mathbf{Y}^N(t),0\le
t<\infty\}$ to $\{\mathbf{y}(t),0\le t<\infty\}$.

{\cor\label{WeakYNtytExamp} In addition to the general assumptions, we assume that $\{\varphi_N\}$, $\varphi$ and
$\{\psi_N\}$, $\psi$ satisfy that for any $T>0$,
\begin{align}
\label{VarphiNconvToVarphiExamp}&\,\lim_{N\to \infty}\sup_{0\le t\le
T}\sup_{\mathbf{x}\in
 K^3_N}|\varphi_N(t,\mathbf{x})-\varphi(t,\mathbf{x})|=0,
 \\
\label{PsiNconvToPsiExamp}&\,\lim_{N\to \infty}\sup_{0\le t\le
T}\sup_{\mathbf{x}\in
 K^3_N}|\psi_N(t,\mathbf{x})-\psi(t,\mathbf{x})|=0.
\end{align}
Suppose that
\begin{equation}\label{eq:ANConvAUniOnFExamp}
\lim_{N\to \infty}\sup_{q\in R}\sup_{\mathbf{x}\in
K^3_N}d_U(A_N(\cdot,\mathbf{x},q), A(\cdot,\mathbf{x}, q))=0
\end{equation}
and that $P(\mathbf{Y}^N(0))^{-1}\Rightarrow\mu$ for some
$\mu\in\mathscr{P}(K^3\times R)$, then there exists on some
probability space $(\Omega, \mathscr{F}, P)$ a stochastic process
$\mathbf{Y}$ satisfying $P(\mathbf{Y}(0))^{-1}=\mu$, which is the
unique solution of the $C_{K^3\times R}[0,\infty)$-martingale
problem for $(G_A, \mu)$ restricted to $C^2_c(K^3\times R)$ and for
$P$-a.s. $\omega\in \Omega$, $\mathbf{Y}(t,\omega)'$ is the unique
solution of (\ref{eq:IAFFMlimitIntegFormx}) and
(\ref{eq:IAFFMlimitIntegFormq}), such that $\mathbf{Y}^N\Rightarrow
\mathbf{Y}$. }

\proof Note that (\ref{VarphiNconvToVarphiExamp}) and
(\ref{PsiNconvToPsiExamp}) imply the Condition
\ref{cond:UniformConvForVarphiPsi}, and
(\ref{eq:ANConvAUniOnFExamp}) implies the Condition
\ref{cond:UniformConvForAN}. \qed

{\rem\label{SufficientCondiForVarphiNVarPsiNPsi} The conditions in
(\ref{VarphiNconvToVarphiExamp}) and (\ref{PsiNconvToPsiExamp}) are
very general. One sufficient condition
to guarantee (\ref{VarphiNconvToVarphiExamp}) and
(\ref{PsiNconvToPsiExamp}) is to assume that for any $T>0$,
$\alpha_i^N(t)$, $\beta_i^N(t)$, $\delta_i^N(t)$ ($1\le i\le 3$), and $F^N(t)$ converges uniformly
 to $\alpha_i(t)$, $\beta_i(t)$,  $\delta_i(t)$, and $F(t)$ on $[0,T]$.
 }

\subsection{Fixed points of the example}\label{sse:FPOASGME}

In this subsection, we consider the fixed point problem for the
example. Make assumptions for
$A(\mathbf{x},q)=(a_{ij}(\mathbf{x},q))_{3 \times 3}$ as follows:
for each $(\mathbf{x},q)\in K^3\times R$,
\begin{itemize}
 \item[1)] $A(\mathbf{x},q) \mathbf{e}={\mathbf 0}_{3 \times 1}$, where $\mathbf{e}=[1, 1, 1]'$;
 \item[2)] $0\le a_{ij}(\mathbf{x},q)\le 1$, for $1\le i,j\le 3$, $i\not=j$;
 \item[3)] $1-\sum_{j\not=i} a_{ij}(\mathbf{x},q)\ge 0$ for $1\le i\le 3$.
\end{itemize}
The above assumptions imply that for each $(\mathbf{x},q)\in K^3\times R$, $E_3+A(\mathbf{x},q)$ is a stochastic
 matrix. Assume also that $A(\mathbf{x},q)$ satisfies the
 Lipschitz condition and is bounded on $K^3\times R$. Then all the conditions related to $A(\mathbf{x},q)$
  required by Lemma \ref{qxqtdxBd} are satisfied, if viewing $A(\mathbf{x},q)$ as a matrix-valued function on
$[0,\infty)\times K^3\times R$.

Next, we assume that in (\ref{varphiExmp}) and (\ref{psiExamp}),
$\alpha_i(t)$, $\beta_i(t)$, and $\delta_i(t)$ ($1\le i\le 3$) and
$\log F(t)$ are constants which are denoted by $\alpha_i$,
$\beta_i$, $\delta_i$ ($1\le i\le 3$) and $\log F$. Assume also that
for any $\mathbf{x}\in K^3$, $x_1 \alpha_1+x_2
\alpha_2(1+\beta_2)+x_3 \alpha_3(1+\beta_3)\not= 0$. Then all the
conditions related to $\varphi(\mathbf{x})$, $\psi(\mathbf{x})$ and
$g(\mathbf{x},q)$ required by Lemma \ref{qxqtdxBd} are satisfied.

To get the fixed points of the following system,
\begin{align}
\label{eq:IAFFMlimitdiffFormxExamp}\frac{d\mathbf{x}(t)}{dt}&=A(\mathbf{x}(t),q(t))' \mathbf{x}(t),  \\
\label{eq:IAFFMlimitdiffFormqExamp}\frac{dq(t)}{dt}&=g(\mathbf{x}(t),q(t)),
\end{align}
we need to solve the following equations:
\begin{align}
\label{eq:FixedPtsAExamp}A(\mathbf{x},q)' \mathbf{x}&=0,  \\
\label{eq:FixedPtsgExamp}g(\mathbf{x},q) &=0.
\end{align}
Note that we use $\mathbf{x}\in K^3$ as a column vector.

At first, we assume that $\alpha_1\not=0$, $\beta_1\not=0$,
$\delta_2\delta_3>0$. We know that $g(\mathbf{x},q)$ is of the
following form:
\begin{equation}\label{eq:gfuncIAFMMexampConstantCoef}
\begin{aligned}
 g(\mathbf{x},q)=&\,\frac{x_1
\alpha_1\beta_1 q}{x_1 \alpha_1+x_2 \alpha_2(1+\beta_2)+x_3
\alpha_3(1+\beta_3)}
\\ &\,+\frac{x_1 \delta_{1}+x_2 \delta_{2}+x_3
\delta_{3}-x_1 \alpha_1\beta_1\log F}{x_1 \alpha_1+x_2
\alpha_2(1+\beta_2)+x_3 \alpha_3(1+\beta_3)}.
\end{aligned}
\end{equation}
Then (\ref{eq:FixedPtsgExamp}) implies that for $\mathbf{x}\in K^3$
with $x_1\not=0$, we have
\begin{equation}\label{eq:Formqx}
q_{\mathbf{x}}=\frac{x_1
\alpha_1\beta_1\log F-x_1 \delta_{1}-x_2 \delta_{2}-x_3
\delta_{3}}{x_1 \alpha_1\beta_1}.
\end{equation}
 Note that for any $\mathbf{x}\in
K^3$ with $x_1=0$, (\ref{eq:FixedPtsgExamp}) implies that
$x_2=x_3=0$ by the assumption $\delta_2\delta_3>0$, which
contradicts with $x_1+x_2+x_3=1$. Therefore, any $\mathbf{x}\in K^3$
with $x_1=0$ does not solve (\ref{eq:FixedPtsgExamp}).

(\ref{eq:Formqx}) defines a function $q_{\mathbf{x}}$:
$K^3\setminus\{\mathbf{x}\in K^3, x_1=0\}\mapsto R$. It is clear
that for any $\tilde{\mathbf{x}}\in\{\mathbf{x}\in K^3, x_1=0\}$,
$\lim_{x_1\not=0,\mathbf{x}\rightarrow \tilde{\mathbf{x}}}
q_{\mathbf{x}}=\infty$ or $-\infty$. We can actually define
$q_{\mathbf{x}}$: $K^3\mapsto R\bigcup\{-\infty,\infty\}$, where for
$\tilde{\mathbf{x}}\in\{\mathbf{x}\in K^3, x_1=0\}$,
 \begin{equation}\label{extendedFuncqx}
 q_{\tilde{\mathbf{x}}}=\lim_{x_1\not=0,\mathbf{x}\rightarrow
\tilde{\mathbf{x}}} q_{\mathbf{x}}.
\end{equation}
 The function $q_{\mathbf{x}}$
with the extended definition is a continuous function on $K^3$.

Since $A(\mathbf{x},q)$ satisfies the Lipschitz condition,
$A(\mathbf{x},q)$ is a continuous function on $K^3\times R$. We make
more assumptions on $A(\mathbf{x},q)$, such that we can extend
$A(\mathbf{x},q)$ to  $K^3\times (R\bigcup \{-\infty, \infty\})$.
Assume that for any $\tilde{\mathbf{x}}\in K^3$,
$\lim_{\mathbf{x}\rightarrow\tilde{\mathbf{x}},q\rightarrow
\pm\infty} A(\mathbf{x},q)$ exists. Then we can extend
$A(\mathbf{x},q)$ to $K^3\times (R\bigcup \{-\infty, \infty\})$ by
defining
\begin{equation}\label{extendedAtoInfinity}
A(\tilde{\mathbf{x}},\pm\infty)=\lim_{\mathbf{x}\rightarrow\tilde{\mathbf{x}},q\rightarrow
\pm\infty} A(\mathbf{x},q),
\end{equation}where $\tilde{\mathbf{x}}\in K^3$. The
extended function $A(\mathbf{x},q)$ is continuous on $K^3\times
(R\bigcup \{-\infty, \infty\})$.

Define for $\mathbf{x}\in K^3$,
\begin{equation}\label{defnTx}
T(\mathbf{x})=A(\mathbf{x}, q_{\mathbf{x}})'\mathbf{x}+\mathbf{x}.
\end{equation}
Since $E_3+A(\mathbf{x}, q_{\mathbf{x}})$ is a stochastic matrix,
$T(\mathbf{x})$ is a map from $K^3$ to $K^3$. Since $q_{\mathbf{x}}$
is a continuous function from $K^3$ to $R\bigcup \{-\infty,
\infty\}$ and $A(\mathbf{x},q)$ is a continuous function from
$K^3\times (R\bigcup \{-\infty, \infty\})$ to $R^{3\times 3}$, we
have that $T(\mathbf{x})$ is a continuous map from $K^3$ to $K^3$.
By Brouwer's fixed point theorem \cite{Smart}, there exists fixed
points $\mathbf{x}^0\in K^3$ such that
$T(\mathbf{x}^0)=\mathbf{x}^0$. It then follows that
\begin{equation}\label{fixedptforA}
A(\mathbf{x}^0, q_{\mathbf{x}^0})'\mathbf{x}^0=0.
\end{equation}
The fixed point of this example might be unique.

We have to make further assumptions on $A(\mathbf{x},q)$ to exclude
the case $x^0_1=0$. Each condition as follows guarantees $x^0_1>0$:
for arbitrary $\tilde{\mathbf{x}}\in\{\mathbf{x}\in K^3, x_1=0\}$,
\begin{enumerate}
\item $a_{21}(\tilde{\mathbf{x}},\pm \infty)>0$ and
      $a_{31}(\tilde{\mathbf{x}},\pm \infty)>0$;
\item $a_{21}(\tilde{\mathbf{x}},\pm \infty)>0$ and
      $a_{32}(\tilde{\mathbf{x}},\pm \infty)>0$;
\item $a_{31}(\tilde{\mathbf{x}},\pm \infty)>0$ and
      $a_{23}(\tilde{\mathbf{x}},\pm \infty)>0$.
\end{enumerate}

The fixed point $\mathbf{x}^0$ determined by (\ref{defnTx}) with
$x^0_1>0$ satisfies that $(\mathbf{x}^0, q_{\mathbf{x}^0})$ is the
solution of the system (\ref{eq:FixedPtsAExamp}) and
(\ref{eq:FixedPtsgExamp}), i.e. $(\mathbf{x}^0, q_{\mathbf{x}^0})$
is the fixed point of (\ref{eq:IAFFMlimitdiffFormxExamp}) and
(\ref{eq:IAFFMlimitdiffFormqExamp}).

{\rem\label{rem:FixedPoint} The set of fixed points is a subset of
the set $\{(\mathbf{x}, q), \mathbf{x}\in K^3, x_1
\alpha_1+\sum_{i=2}^3 x_i \alpha_i(1+\beta_i)= 0, x_1
\alpha_1\beta_1(q-\log F) +\sum_{i=1}^3 x_i \delta_{i}=0\}$. The
latter set contains the stationary solutions or explosive solutions
of our financial system. }

\subsection{Connection with classical stock price formula}\label{sse:CWCSPF}

We make the same assumptions on $\alpha_i^N(t)$, $\beta_i^N(t)$, $\alpha_i(t)$ and $\beta_i(t)$ ($1\le i\le 3$)
 as those in subsection \ref{sse:CAVOTGME}. But we assume that $\delta_i^N(t)$, ($1\le i\le 3$), $F^N(t)$
 defined on some probability space $(\Omega_N, \mathscr{F}_N, P_N)$, with sample pathes satisfying the conditions
 in subsection \ref{sse:CAVOTGME} $P_N$-a.s.; and that $\delta_i(t)$ ($1\le i\le 3$) are Brownian motions,  and $F(t)$ is a geometric Brownian motion.
 Assume that $(\delta^N_1, \delta^N_2, \delta^N_3, F^N)\Rightarrow (\delta_1,
\delta_2, \delta_3, F)$ in the sense of weak convergence. Assume
also all other conditions for $A_N(t, \mathbf{x}, q)$ and $A(t,
\mathbf{x}, q)$ used in subsection \ref{sse:CAVOTGME}.

{\rem\label{rem:liqu&fundv} Duffie and Protter (see
\cite{Duffie&Protter}) justified the weak convergence of properly
scaled liquidity demand to Brownian motions.
 It is also usual to assume the fundamental value of a financial asset to be a geometric Brownian motion.
 The liquidity demands $\delta_i(t)$ ($1\le i\le 3$) and the fundamental value $F(t)$ constitute the random environment for the limit process $Y$. }

We can follow the procedures in subsection \ref{sse:AGMEOIAFFM} and
\ref{sse:CAVOTGME} to define $\{\mathbf{Y}^N(t,\omega_N),\newline
0\le t<\infty\}$ and determine $\{\mathbf{Y}(t,\omega), 0\le
t<\infty\}$. As to the determination of $\{\mathbf{Y}(t,\omega),
0\le t<\infty\}$ by the {\em nonlinear Volterra Equations of the
second kind}, we note that the Brownian paths $\delta_i(t,\omega)$
($1\le i\le 3$) and $\log F(t,\omega)$ satisfy
(\ref{BdExpGrowthDeltaF}) almost surely by the Law of the iterated
logarithm. Then we can prove that $(\delta_1^N,\delta_2^N,\delta_3^N,F^N, \mathbf{Y}^N)\Rightarrow
(\delta_1,\delta_2,\delta_3,F, \mathbf{Y})$, see the method used in \cite{Wu2a}. As usual, the Skorohod
Representation theorem is the basic tool to establish this result.
Therefore, without loss of generality, we can assume that
$(\delta^N_1, \delta^N_2, \delta^N_3, F^N)$ and $(\delta_1,
\delta_2, \delta_3, F)$ are defined in the same probability space
$(\Omega, \mathscr{F}, P)$ and $(\delta^N_1, \delta^N_2, \delta^N_3,
F^N)$ converges to $(\delta_1, \delta_2, \delta_3, F)$ almost
surely. Note that for $1\le i\le 3$, $\lim_{N\to
\infty}d(\delta_i^N(\omega),\delta_i(\omega))=0$ implies that
$\lim_{N\to \infty}d_U(\delta_i^N(\omega),\delta_i(\omega))=0$ and
hence $\delta_i^N(\omega)$ converges to $\delta_i(\omega)$ uniformly
on $[0, T]$ for any $T>0$. Same thing holds for $F^N(\omega)$ and
$F(\omega)$. Then we actually assumed or verified all the conditions
required by Remark \ref{SufficientCondiForVarphiNVarPsiNPsi}.

Note that the log price function $q(t)$ in random environment is
determined pathwisely by the {\em nonlinear Volterra Equations of
the second kind}. That is say, $P$-a.s for any sample point
$\omega\in \Omega$, the {\em nonlinear Volterra Equations of the
second kind} determines a unique log price function $q(t, \omega)$
on $[0, \infty)$. This is similar to the classical assumption of the
stock price formula which was suggested by Samuelson
\cite{Samuelson} in 1964: $\log S_t =\mu t +\sigma W_t$, where $W_t$
is a Brownian motion. Therefore, the interacting agent feedback
financial system in random environment
 generalizes the classical stock price formula by
incorporating the interaction between different types of agents and
the interaction between the stock price and the empirical
distribution of the types of agents.

\section{Proof of Theorem \ref{thm:ConvIAFFM}}\label{se:ProofofThm}

At first, we consider the moments for $\{ \mathbf{n}^N(k), k \ge
0\}$. Let $\mathbf{V}=(v_1, \cdots, v_r)'$ be a positive
 vector. For any $k\ge 1$, by (\ref{eq:IAFFMXiNkj}), (\ref{IAFFMeq:nk+In}), and the independence of $\mathbf{\Xi}_{N,k,i}$'s, we have
\begin{equation}\label{eq:mgfforIAFFM}
E\biggl[\prod_{i=1}^{r}v_i^{n^N_i(k+1)}\biggl|(\mathbf{n}^N(k),\tilde{q}^N(\frac{k}{N}))\biggl]=\prod_{i=1}^{r}\biggl(P_{N,i\cdot}[k,\mathbf{X}^N(\frac{k}{N}),\tilde{q}^N(\frac{k}{N})]\mathbf{V}\biggl)^{n^N_i(k)}.
\end{equation}
It follows by (\ref{eq:mgfforIAFFM}) that
\begin{equation}\label{eq:meanIAFFM}
E\biggl[{\mathbf n}^N(k+1)\biggl|({\mathbf
n}^N(k),\tilde{q}^N(\frac{k}{N}))\biggl]={\mathbf
n}^N(k)P_{N}[k,\mathbf{X}^N(\frac{k}{N}),\tilde{q}^N(\frac{k}{N})],
\end{equation}
and for $1\le i\le r$,
\begin{equation}\label{eq:squMinusMeanIAFFM}
\begin{aligned}
&\,E\biggl[n^N_i(k+1)(n^N_i(k+1)-1)\biggl|({\mathbf n}^N(k),\tilde{q}^N(\frac{k}{N}))\biggl] \\
=&\,\biggl({\mathbf n}^N(k)P_{N,\cdot,i}[k,\mathbf{X}^N(\frac{k}{N}),\tilde{q}^N(\frac{k}{N})]\biggl)^2 \\
  &\, -\sum_{j=1}^{r}\biggl(P_{N,j\cdot}[k,\mathbf{X}^N(\frac{k}{N}),\tilde{q}^N(\frac{k}{N})]P_{N,\cdot,i}[k,\mathbf{X}^N(\frac{k}{N}),\tilde{q}^N(\frac{k}{N})]\biggl)^2 n^N_{j}(k),
\end{aligned}
\end{equation}
 where
$P_{N,\cdot,i}[k,\mathbf{X}^N(\frac{k}{N}),\tilde{q}^N(\frac{k}{N})]$
is $i$-th column of
$P_{N}[k,\mathbf{X}^N(\frac{k}{N}),\tilde{q}^N(\frac{k}{N})]$.
 Then we can get
\begin{equation}\label{eq:squareOfDifIAFFM}
\begin{aligned}
&\,E[(n^N_i(k+1)-n^N_i(k))^2\biggl|({\mathbf n}^N(k),\tilde{q}^N(\frac{k}{N}))\biggl]\\
=&\,\biggl({\mathbf n}^N(k)P_{N,\cdot,i}[k,\mathbf{X}^N(\frac{k}{N}),\tilde{q}^N(\frac{k}{N})]\biggl)^2 +{\mathbf n}^N(k)P_{N,\cdot,i}[k,\mathbf{X}^N(\frac{k}{N}),\tilde{q}^N(\frac{k}{N})] \\
&\,
-\sum_{j=1}^{r}\biggl(P_{N,j\cdot}[k,\mathbf{X}^N(\frac{k}{N}),\tilde{q}^N(\frac{k}{N})]P_{N,\cdot,i}[k,\mathbf{X}^N(\frac{k}{N}),\tilde{q}^N(\frac{k}{N})]\biggl)^2
n^N_{j}(k)  \\
&\,-2{\mathbf
n}^N(k)P_{N,\cdot,i}[k,\mathbf{X}^N(\frac{k}{N}),\tilde{q}^N(\frac{k}{N})]n^N_i(k)+(n^N_i(k))^2.
\end{aligned}
\end{equation}

Next, we prove that $\{\mathbf{Y}^N\}$ satisfies the compact
containment condition under certain conditions on
$\{\mathbf{Y}^N(0)\}$, $\{g_N\}$ and $g$.

{\lem\label{lem:CompactContainForYN} Assume that
$P(\mathbf{Y}^N(0))^{-1}\Rightarrow\mu$ for some
$\mu\in\mathscr{P}(K\times R)$, and  $g(t,\mathbf{x},q)$,
$g_N(t,\mathbf{x},q)$ satisfy (\ref{eq:FormOfg}) and
(\ref{eq:FormOfgN}). Assume also that for any $T>0$,
$\varphi(t,\mathbf{x}, q)$ and $\psi(t,\mathbf{x}, q)$ are bounded
on $[0,T]\times K\times R$ and the condition
\ref{cond:UniformConvForVarphiPsi} holds. Then $\{\mathbf{Y}^N\}$
satisfies the compact containment condition, i.e., for every
$\eta>0$ and $T>0$, there exists a compact set
$\tilde{K}_{\eta,T}\subset K\times R$ for which
\begin{equation}\label{eq:CompactContIAFFM}\inf_{N}P\{\mathbf{Y}^N(t)\in
\tilde{K}_{\eta,T},\mbox{ for } 0\le t\le T\}\ge 1-\eta.
\end{equation}

\proof Since  $P(\mathbf{Y}^N(0))^{-1}\Rightarrow\mu$,
$\{P(q^N(0))^{-1}\}$ is tight. For any $\eta>0$, there exists $b>0$
such that
\begin{equation*}
 \inf_N P\{q^N(0)\in I\}\ge 1-\eta,
\end{equation*}
 where $I=[-b,b]$.

Fix $T>0$, let $C_T>0$ be a bound of $\varphi(t,\mathbf{x},q)$ and
$\psi(t,\mathbf{x},q)$ on $[0,T]\times K\times R$. By
(\ref{eq:RecLogPriceTdqNgN}), (\ref{eq:FormOfg}),
(\ref{eq:FormOfgN}), (\ref{VarphiNconvToVarphi}) and
(\ref{PsiNconvToPsi}), there exists $N_0$, such that for $N>N_0$, we
have for any $0\le k\le [NT]-1$,
\begin{equation}\label{eq:bdqNk1minusqNk}
|\tilde{q}^N(\frac{k+1}{N})-\tilde{q}^N(\frac{k}{N})|\le
\frac{1}{N}C_T[|\tilde{q}^N(\frac{k}{N})|+1],
\end{equation}
 which implies that for
$0\le k\le [NT]$,
\begin{equation}\label{eq:bdqNk}
|\tilde{q}^N(\frac{k}{N})|\le (1+\frac{1}{N}C_T)^k
[|\tilde{q}^N(0)|+1]-1.
\end{equation}
 Let $I_{\eta,T}=[-e^{T
C_T}(b+1), e^{T C_T}(b+1)]$. Since $(1+\frac{1}{N}C_T)^N$ increases
as $N$ does, with the limit $e^{C_T}$, we have for $N> N_0$
\begin{equation*}
 P\{q^N(t)\in I_{\eta,T}, 0\le t\le T\}=P\{q^N(0)\in I\}.
\end{equation*}

Let $\tilde{K}_{\eta, T}=K\times I_{\eta,T}$, then for $N> N_0$
\begin{equation*}
P\{\mathbf{Y}^N(t)\in \tilde{K}_{\eta,T}\}=P\{q^N(0)\in I\}
\end{equation*}
 and
(\ref{eq:CompactContIAFFM}) holds. \qed

Next, we state Lemma \ref{rem:verifyCondiEofCor3.6Chp2} which
verifies (3.10) in condition (\textbf{e}) of Corollary 3.5 for
$\{\mathbf{Y}^N\}$, see \cite{Wu2b}. The proof of Lemma
\ref{rem:verifyCondiEofCor3.6Chp2} is given in Appendix A.

{\lem\label{rem:verifyCondiEofCor3.6Chp2} Assume that $A(t,
\mathbf{x}, q)=(a_{i,j}(t, \mathbf{x}, q))_{r\times r}$ satisfies
the condition \ref{cond:regA} and
$\{A_N(t,\mathbf{x},q)=(a_{N,i,j}(t, \mathbf{x}, q))_{r\times r}\}$
satisfies the condition \ref{cond:regAN}. Assume that $A(t,
\mathbf{x}, q)$ is bounded on $[0,T]\times K\times R$ for any $T>0$.
Assume that $g(t,\mathbf{x},q)$ and $\{g_N(t,\mathbf{x},q)\}$
satisfy (\ref{eq:FormOfg}) and (\ref{eq:FormOfgN}). Assume that
$\varphi$, $\psi$, $\frac{\partial \varphi}{\partial {\mathbf x}}$,
$\frac{\partial \psi}{\partial {\mathbf x}}$ are continuous, and the
condition \ref{cond:UniformConvForVarphiPsi} holds for $\{
\varphi_N\}$, $\varphi$ and $\{\psi_N\}$, $\psi$. If for each $f\in
C^2_c(K\times R)$ and $T>0$ there exist measurable sets
$\{F_N\}\subset R$ such that (\ref{eq:ANConvAUniOnF}) holds, then
\begin{equation}\label{GNfGfIAFFM}
\lim_{N\to \infty}\sup_{0\le t\le T}\sup_{q\in
F_N}\sup_{\mathbf{x}\in
K_N}|N[S_{N,[Nt]}-I]f(\mathbf{x},q)-G_A(t)f(\mathbf{x},q)|=0.
 \end{equation}
}

{\rem\label{specialgN} As a special case, we can assume that $g_N(t,
\mathbf{x},q)=g(\frac{[Nt]}{N}, \mathbf{x},q)$, i.e., $\varphi_N(t,
\mathbf{x},q)=\varphi(\frac{[Nt]}{N}, \mathbf{x},q)$ and $\psi_N(t,
\mathbf{x},q)=\psi(\frac{[Nt]}{N}, \mathbf{x},q)$, for $(t,
\mathbf{x},q)\in [0, \infty)\times K_N\times R$. Then we do not need
to assume the Condition \ref{cond:UniformConvForVarphiPsi} to get
the conclusions of Lemma \ref{lem:CompactContainForYN} and Lemma
\ref{rem:verifyCondiEofCor3.6Chp2}. But the Condition
\ref{cond:UniformConvForVarphiPsi} holds under the stronger
condition that $\frac{\partial \varphi}{\partial t}$ and
$\frac{\partial \psi}{\partial t}$ are continuous and bounded on
$[0, T]\times K\times R$ for any $T>0$. }

{\cor\label{cor:YNtight} If we assume in Lemma
\ref{rem:verifyCondiEofCor3.6Chp2} that
$P(\mathbf{Y}^N(0))^{-1}\Rightarrow\mu$ for some
$\mu\in\mathscr{P}(K\times R)$, $g(t,\mathbf{x},q)$ and
$\{g_N(t,\mathbf{x},q)\}$ satisfy (\ref{eq:FormOfg}) and
(\ref{eq:FormOfgN}), $\varphi(t,\mathbf{x}, q)$ and
$\psi(t,\mathbf{x}, q)$ are bounded on $[0,T]\times K\times R$ for
any $T>0$, and the Condition \ref{cond:UniformConvForAN} holds, then
$\{\mathbf{Y}^N\}$ is relatively compact. If we assume also that the
Condition \ref{cond:ExpGrowth} holds, then any limit point
$\mathbf{Y}$ of $\mathbf{Y}^N$ is a solution of the $D_{K\times
R}[0,\infty)$-martingale problem for $(G_A, \mu)$ restricted to
$C^2_c(K\times R)$. }


\proof By Lemma \ref{lem:CompactContainForYN}, we know that
$\{\mathbf{Y}^N\}$ satisfies the compact containment condition. Let
$\tilde{F}_N=K_N\times F_N$, by (\ref{eq:YnFnConInProbOneIFMM}) and
Lemma \ref{rem:verifyCondiEofCor3.6Chp2}, we know that (3.9) and (3.10)
in the condition (\textbf{e}) of Corollary 3.5 \cite{Wu2b}, hold for $\{\tilde{F}_N\}$. Note that
 \begin{equation*}
 \|G_A(t)f\|\le \|\mathbf{x}A(t,\mathbf{x},q)\frac{\partial f}{\partial
\mathbf{x}}\|+\|\frac{\partial f}{\partial q}g(t,\mathbf{x},q)\|,
\end{equation*}
where $\|\cdot\|$ above is the supnorm with respect to
$(\mathbf{x},q)$, and that $A(t, \mathbf{x}, q)$ is bounded on
$[0,T]\times K\times R$ for any $T>0$, $g$ is continuous on $[0,
\infty)\times K\times R$, and $f$ has a compact support, we have
that $\|G_A(t)f\|$ is bounded on $[0, T]$ for any $T>0$. Then by the
proof of Corollary 3.5 \cite{Wu2b}, we
know that the equations (2.5) and
(2.6) of \cite{Wu2b} are verified.
Then it follows that $\{\mathbf{Y}^N\}$ is relatively compact.

Let $\tilde{K}$ be the support of $f$, then by (\ref{eq:FormOfg}) we
get
\begin{equation}\label{eq:GAtfbdExp}
  \begin{aligned}
 \|G_A(t)f\| \le & \,\max_{1\le i\le r}\sup_{(\mathbf{x},q)\in
\tilde{K}}\|\frac{\partial f}{\partial x_i}\|
\cdot\sup_{(\mathbf{x},q)\in
 \tilde{K}}\|A(t,\mathbf{x},q)\|
 +\sup_{(\mathbf{x},q)\in \tilde{K}}|\frac{\partial
f}{\partial q}q|\cdot\sup_{(\mathbf{x},q)\in
\tilde{K}}|\varphi(t,\mathbf{x},q)| \\
 &\,+\sup_{(\mathbf{x},q)\in
\tilde{K}}|\frac{\partial f}{\partial
q}|\cdot\sup_{(\mathbf{x},q)\in \tilde{K}}|\psi(t,\mathbf{x},q)|,
\end{aligned}
\end{equation}
where $\|G_A(t)f\|$ is the supnorm of $G_A(t)f$ with respect to
$(\mathbf{x},q)$, and the norms on the right side of
(\ref{eq:GAtfbdExp}) are matrix norm and Euclidean norm used in the
proof of  Lemma \ref{rem:verifyCondiEofCor3.6Chp2}. By
(\ref{eq:IAFFMA(s)ExpGrowth}), (\ref{eq:IAFFMVarphi(s)ExpGrowth})
and (\ref{eq:IAFFMPsi(s)ExpGrowth}), we know that (2.7) in the
condition (\textbf{b}) of Proposition 2.1 \cite{Wu2b} is verified.
Then by the proof of Corollary 3.5 \cite{Wu2b}, any limit point
$\mathbf{Y}$ of $\mathbf{Y}^N$ is a solution of the $D_{K\times
R}[0,\infty)$-martingale problem for $(G_A, \mu)$ restricted to
$C^2_c(K\times R)$. \qed

Now we state Lemma \ref{lem:LimitPtBeingCont} and give its proof in Appendix B.
{\lem\label{lem:LimitPtBeingCont} If we assume all the assumptions
made in Lemma \ref{rem:verifyCondiEofCor3.6Chp2} and Corollary
\ref{cor:YNtight}, then any limit point $Y$ of $\mathbf{Y}^N$ is a
solution of the $C_{K\times R}[0,\infty)$-martingale problem for
$(G_A, \mu)$ restricted to $C^2_c(K\times R)$. }

It is clear that any solution of the $C_{K\times
R}[0,\infty)$-martingale problem for $(G_A, \mu)$ restricted to
$C^2_c(K\times R)$ is a solution of the integral equations
(\ref{eq:IAFFMlimitIntegFormx}) and (\ref{eq:IAFFMlimitIntegFormq}).

{\rem\label{rem:LipschitzBigb} If $A(t,\mathbf{y})$ satisfies the
Lipschitz condition, then
$\hat{\mathbf{b}}(t,\mathbf{y})=\mathbf{x}A(t,\mathbf{y})$ satisfies
the Lipschitz condition under the condition that $A(t,\mathbf{y})$
is bounded on $[0,T]\times K\times R$ for any $T>0$. But even if
$\varphi(t,\mathbf{y})$ and $\psi(t,\mathbf{y})$ are assumed to be
bounded on $[0,T]\times K\times R$ for any $T>0$ and satisfy the
Lipschitz condition, we can not claim that
$b_{r+1}(t,\mathbf{y})=g(t,\mathbf{y})$ satisfies the Lipschitz
condition. This is the reason that we need to expand the standard
theory of {\em nonlinear Volterra equations of the second kind} for
$g(t,\mathbf{y})$ that does not satisfy the Lipschitz condition.}

We need to introduce notations and property related to the {\em
nonlinear Volterra equations of the second kind} with semi-Lipschitz
conditions. Let $V(s,t,\mathbf{y})$ be a measurable $R^m$-valued
function on $[a,b]\times[a,b]\times R^m$ and $f$ an
$\mathbb{L}^2([a,b], R^{m})$ function. Let $m=n+1$. We give
different notations for the first $n$ components of $\mathbf{y}$,
$V(s,t,\mathbf{y})$ and $f(s)$. For any vector
$\mathbf{y}=(y_1,\cdots,y_n,y_{n+1})'\in R^{n+1}$, we let
$\mathbf{x}=(y_1,\cdots,y_n)'$. Similarly, we let
$\hat{V}(s,t,\mathbf{y})=(V_1(s,t,\mathbf{y}),\cdots,V_n(s,t,\mathbf{y}))'$
and $\hat{f}(s)=(f_1(s),\cdots,f_n(s))'$. Then we can represent the
integral equations with kernel $V(s,t,\mathbf{y})$ and a function
$f$ as follows: for $a\le s\le b$,
\begin{align}
\label{NonlinearVoltEqFirstn}\mathbf{x}(s)&=\hat{f}(s)+\int_a^s \hat{V}(s,t,\mathbf{x}(t),y_{n+1}(t))dt  \\
\label{NonlinearVoltEqLast}y_{n+1}(s)&=f_{n+1}(s)+\int_a^s
V_{n+1}(s,t,\mathbf{x}(t),y_{n+1}(t))dt
\end{align}

The semi-Lipschitz conditions on (\ref{NonlinearVoltEqFirstn}) and
(\ref{NonlinearVoltEqLast}) are as follow:
\begin{itemize}
\item[\textbf{(I)}] If $(\tilde{\mathbf{x}}(t)',\tilde{y}_{m+1}(t))'$ is any
solution of (\ref{NonlinearVoltEqFirstn}) and
(\ref{NonlinearVoltEqLast}), then $\tilde{y}_{m+1}(t)$ is the unique
solution of (\ref{NonlinearVoltEqLast}) if we plug into
(\ref{NonlinearVoltEqLast}) $\mathbf{x}(t)=\tilde{\mathbf{x}}(t)$.
That is to say, there exists a dependence, which is based on
(\ref{NonlinearVoltEqLast}), between the first $n$ components and
the last component of any solution. Then we can actually write
$\tilde{y}_{m+1}(t)=\xi_{\tilde{\mathbf{x}}}(t)$.
\item[\textbf{(II)}] $\hat{V}(s,t,\mathbf{y})$ satisfies the Lipschitz condition, i.e. there exists a
nonnegative $\mathbb{L}^2([a,b]\times [a,b], R)$ function
$V_0(s,t)$, such that
 for any $a\le t<s\le b$ and $\mathbf{y}_1=(\mathbf{x}_1',y_{n+1}^{(1)})',
\mathbf{y}_2=(\mathbf{x}_2',y_{n+1}^{(2)})' \in R^{n+1}$,
\begin{equation}\label{HatVLipschitz}
\|\hat{V}(s,t,\mathbf{y}_1)-\hat{V}(s,t,\mathbf{y}_2)\|\le V_0(s,t)
\sqrt{\|\mathbf{x}_1-\mathbf{x}_2\|^2+(y_{n+1}^{(1)}-y_{n+1}^{(2)})^2}.
\end{equation}
\item[\textbf{(III)}] Based on \textbf{(I)}, if $(\mathbf{x}(t)',\xi_{\mathbf{x}}(t))'$ and
 $(\tilde{\mathbf{x}}(t)',\xi_{\tilde{\mathbf{x}}}(t))'$ are two
solutions of (\ref{NonlinearVoltEqFirstn}) and
(\ref{NonlinearVoltEqLast}), then there exists a nonnegative
measurable function $B(s)$ on $[a,b]$ such that
 \begin{equation}\label{DiffOfTwoSolns}
 \int_a^s (\xi_{\mathbf{x}}(t)-\xi_{\tilde{\mathbf{x}}}(t))^2 dt\le B(s)
 \int_a^s\|\mathbf{x}(t)-\tilde{\mathbf{x}}(t)\|^2 dt \mbox{ for }
 a\le s\le b,
 \end{equation}
and
\begin{equation}\label{V0stBtInteg}
\int_a^b [\int_a^s V_0(s,t)^2 dt (1+B(s))]ds
 <\infty.
 \end{equation}
\end{itemize}

{\prop\label{NonlinearVoltSemiLips} Assume that $V$ is
 a measurable $R^{n+1}$-valued function on $[a,b]\times[a,b]\times R^{n+1}$
 satisfying the semi-Lipschitz conditions \textbf{(I)}-\textbf{(III)} above. Then the nonlinear
Volterra equations of the second kind (\ref{NonlinearVoltEqFirstn})
and (\ref{NonlinearVoltEqLast}) have at most one
$\mathbb{L}^2([a,b], R^{n+1})$ solution. }

\proof Assume that $V_0(s,t)$ is a nonnegative
$\mathbb{L}^2([a,b]\times [a,b], R)$ function which satisfies
(\ref{HatVLipschitz}). Assume that
$(\mathbf{x}(t)',\xi_{\mathbf{x}}(t))'$ and
 $(\tilde{\mathbf{x}}(t)',\xi_{\tilde{\mathbf{x}}}(t))'$ are two
solutions of (\ref{NonlinearVoltEqFirstn}) and
(\ref{NonlinearVoltEqLast}) and $B(s)$ is a nonnegative measurable
function on $[a,b]$ which satisfies (\ref{DiffOfTwoSolns}) and
(\ref{V0stBtInteg}). Let $\mathscr{A}(s)^2=\int_a^s V_0(s,t)^2 dt
(1+B(s))$ for $a\le s\le b$ and let $h^2=\mathscr{A}(b)^2$. Then by
(\ref{HatVLipschitz}), Cauchy-Schwartz inequality and
(\ref{V0stBtInteg}), we have for $a\le s\le b$,
\begin{equation}\label{eq:diffXsTdXs}
\begin{aligned}\|\mathbf{x}(s)-\tilde{\mathbf{x}}(s)\|^2\le
&\,\{\int_a^s\|\hat{V}(s,t,\mathbf{x}(t),\xi_{\mathbf{x}}(t))
-\hat{V}(s,t,\tilde{\mathbf{x}}(t),\xi_{\tilde{\mathbf{x}}}(t))\|dt\}^2
\\
\le &\,\int_a^s V_0^2(s,t)dt \int_a^s
[\|\mathbf{x}(t)-\tilde{\mathbf{x}}(t)\|^2+
 (\xi_{\mathbf{x}}(t)-\xi_{\tilde{\mathbf{x}}}(t))^2 ]dt \\
 \le &\,\int_a^s V_0^2(s,t)dt \biggl[\int_a^s
\|\mathbf{x}(t)-\tilde{\mathbf{x}}(t)\|^2 dt+ B(s)\int_a^s
\|\mathbf{x}(t)-\tilde{\mathbf{x}}(t)\|^2 dt\biggl] \\
 =&\, \mathscr{A}(s)^2 \int_a^s
\|\mathbf{x}(t)-\tilde{\mathbf{x}}(t)\|^2 dt.
 \end{aligned}
\end{equation}

Put $k^2=\int_a^b \|\mathbf{x}(t)-\tilde{\mathbf{x}}(t)\|^2 dt$. By
successive substitutions into (\ref{eq:diffXsTdXs}), we can get
\begin{equation}\label{eq:IntegdiffXsTdXs}
\int_a^s\|\mathbf{x}(t)-\tilde{\mathbf{x}}(t)\|^2 dt\le k^2
\frac{1}{l!}\biggl[\int_a^s  \mathscr{A}(t)^2 dt\biggl]^l\le k^2
\frac{h^{2l}}{l!}
\end{equation}
 for any integer $l\ge 1$ and $a\le s\le b$. Let
$l\rightarrow \infty$, we conclude that $\mathbf{x}(s)\equiv
\tilde{\mathbf{x}}(s)$ for $s\in [a,b]$ in the sense of
$\mathbb{L}^2([a,b], R^{n})$. By condition \textbf{(I)}, we also get
that $\xi_{\mathbf{x}}(s)\equiv\xi_{\tilde{\mathbf{x}}}(s)$ for
$s\in [a,b]$ in the sense of $\mathbb{L}^2([a,b], R)$. \qed

{\lem\label{lem:UniqForIAFFM} Assume that $A(t, \mathbf{x},
q)=(a_{i,j}(t, \mathbf{x}, q))_{r\times r}$ satisfies the Condition
\ref{thm:ConvIAFFM} and is bounded on $[0,T]\times K\times R$ for
any $T>0$. Assume that $g(t,\mathbf{x},q)$ satisfy
(\ref{eq:FormOfg}). Assume also that $\varphi(t, \mathbf{x}, q)$ and
$\psi(t, \mathbf{x}, q)$ are bounded on $[0,T]\times K\times R$ for
any $T>0$ and that $\varphi$, $\psi$, $\frac{\partial
\varphi}{\partial {\mathbf x}}$, $\frac{\partial \psi}{\partial
{\mathbf x}}$ are continuous. Define $\mathbf{b}(t,\mathbf{y})$ and
$\hat{\mathbf{b}}(t,\mathbf{y})$ by (\ref{defnbiForOneTor}),
(\ref{defnbr+1}) and (\ref{eq:BigBtHatBt}). Assume that either
$\mathbf{b}(t,\mathbf{y})$ satisfies the Lipschitz condition or that
the condition \ref{cond:semiLipForIAFFM} holds for
(\ref{eq:IAFFMlimitIntegFormx}) and (\ref{eq:IAFFMlimitIntegFormq}).
 Then the $C_{K\times R}[0,\infty)$-martingale problem for $(G_A, \mu)$ restricted
 to $C^2_c(K\times R)$ has at most one solution. }

\proof First, we prove the case when $\mathbf{b}(t,\mathbf{y})$
satisfies the Lipschitz condition. We just need to prove that the
deterministic integral equation (\ref{eq:DetermInteEq}) has at most
one solution. Define
 \begin{equation}\label{DefnOfVMFFforUniq}
 V(s,t,\mathbf{y})= \left\{
\begin{array}{ll}
                           \mathbf{b}(t,\mathbf{y})' & \mbox{if $0\le t\le s<T$}, \\
                                 0,         & \mbox{if $s,t\in [0, T]$ and $t>s$.}
 \end{array} \right.
 \end{equation}
Note that any solution of the integral equation
(\ref{eq:DetermInteEq}) is continuous. By the Lipschitz condition on
$\mathbf{b}(t,\mathbf{y})$, the {\em nonlinear Volterra equations of
the second kind} (\ref{eq:DetermInteEq}) has only one continuous
solution, see \cite{Tricomi}.

Second, we prove the case when (\ref{eq:IAFFMlimitIntegFormx}) and
(\ref{eq:IAFFMlimitIntegFormq}) satisfy the semi-Lipschitz
condition. This proof is similar to the first case. We just need to
prove that the deterministic integral equations
(\ref{eq:IAFFMlimitIntegFormx}) and (\ref{eq:IAFFMlimitIntegFormq})
have at most one solution. Define
 \begin{equation}\label{DefnOfHatVMFFforUniq}
 \hat{V}(s,t,\mathbf{y})= \left\{
\begin{array}{ll}
                           \hat{\mathbf{b}}(t,\mathbf{y})' & \mbox{if $0\le t\le s<T$}, \\
                                 0,         & \mbox{if $s,t\in [0, T]$ and $t>s$;}
 \end{array} \right.
 \end{equation}
and
 \begin{equation*}
 V_{r+1}(s,t,\mathbf{y})=\left\{
\begin{array}{ll}
 g(t,\mathbf{y}) & \mbox{if $0\le t\le s<T$,} \\
  0,         & \mbox{if $s,t\in [0, T]$ and $t>s$.}
 \end{array} \right.
 \end{equation*}
Note that any solution of the integral equations
(\ref{eq:IAFFMlimitIntegFormx}) and (\ref{eq:IAFFMlimitIntegFormq})
is continuous. We can check easily that the semi-Lipschitz
conditions \textbf{(I)}-\textbf{(III)} for $V(s,t,\mathbf{y})$
required by Proposition \ref{NonlinearVoltSemiLips} are satisfied.
Then (\ref{eq:IAFFMlimitIntegFormx}) and
(\ref{eq:IAFFMlimitIntegFormq}) have at most one continuous
solution.          \qed

{\rem\label{RemkOnFucg(t,y)} The requirement on $g$ that for fixed
$T>0$ and any given continuous $\mathbf{x}(t)$ on $[0,T]$,
(\ref{eq:IAFFMlimitIntegFormq}) determines a unique solution
$q_{\mathbf{x}}(t)$ on $[0,T]$ looks like a awkward one, but it can
be easily satisfied for certain kind of functions. For example, $g$
satisfies this requirement if we assume that the variable $q$ in
$\varphi(t, \mathbf{x}, q)$ and $\psi(t, \mathbf{x}, q)$ is a dummy
variable. }

{\textbf{Proof of Theorem \ref{thm:ConvIAFFM}} } \hspace{2mm} The
conclusion follows by Corollary \ref{cor:YNtight}, Lemma
\ref{lem:LimitPtBeingCont}, and Lemma \ref{lem:UniqForIAFFM}.  \qed


\section*{Appendix A: Proof of Lemma \ref{qxqtdxBd}}

For fixed $T>0$ and any given continuous $K^3$-valued function
$\mathbf{x}(t)$ on $[0,T]$, we define $P_{\mathbf{x}}(t)$ and
$Q_{\mathbf{x}}(t)$ as follows:
\begin{align}
\label{PxtExmp}&\,P_{\mathbf{x}}(t)=\frac{x_1(t)
\alpha_1(t)\beta_1(t)}{x_1(t) \alpha_1(t)+x_2(t)
\alpha_2(t)(1+\beta_2(t))+x_3(t) \alpha_3(t)(1+\beta_3(t))},  \\
\label{QxtExamp} &\, Q_{\mathbf{x}}(t)=\frac{x_1(t)
\delta_{1}(t)+x_2(t) \delta_{2}(t)+x_3(t) \delta_{3}(t)-x_1(t)
\alpha_1(t)\beta_1(t)\log F(t)}{x_1(t) \alpha_1(t)+x_2(t)
\alpha_2(t)(1+\beta_2(t))+x_3(t) \alpha_3(t)(1+\beta_3(t))}.
\end{align}
Then (\ref{dqtdtgtxq}) becomes
\begin{equation}\label{dqtdtgtxqExamp}
\frac{dq(t)}{dt}=P_{\mathbf{x}}(t) q+Q_{\mathbf{x}}(t).
\end{equation}
The unique solution of (\ref{dqtdtgtxqExamp}) with initial condition
$q(0)=q_{\mathbf{x}}(0)$ is given by
\begin{equation}\label{FormOfqxt}
q_{\mathbf{x}}(t)=e^{\int_0^t
P_{\mathbf{x}}(u)du}\biggl[\int_0^t Q_{\mathbf{x}}(u)e^{-\int_0^u
P_{\mathbf{x}}(v)dv}du+q_{\mathbf{x}}(0)\biggl].
\end{equation}
Let $h_{\mathbf{x}}(t)=x_1(t) \alpha_1(t)+x_2(t)
\alpha_2(t)(1+\beta_2(t))+x_3(t) \alpha_3(t)(1+\beta_3(t))$, the
denominator of $P_{\mathbf{x}}(t)$. Since $\alpha_i(t)$,
$\beta_i(t)$ ($1\le i\le 3$) are continuous and for any
$\mathbf{x}\in K^3$ and $t\ge 0$, $x_1 \alpha_1(t)+x_2
\alpha_2(t)(1+\beta_2(t))+x_3 \alpha_3(t)(1+\beta_3(t))\not= 0$, it
follows that
\begin{equation}\label{bdhxt}
B_L\le |h_{\mathbf{x}}(t)|\le B_U, \mbox{ for all $t\in [0,
T]$,}
\end{equation}
where $B_L$, $B_U>0$ do not depend on $\mathbf{x}(t)$.

We assume that $(\mathbf{x}(t),q_{\mathbf{x}}(t))$ and
$(\tilde{\mathbf{x}}(t), q_{\tilde{\mathbf{x}}}(t))$ are two
solutions of (\ref{eq:IAFFMlimitIntegFormx}) and (\ref{eq:IAFFMlimitIntegFormq}) on $[0,T]$. It follows that
$\mathbf{x}(0)=\tilde{\mathbf{x}}(0)$ and
$q_{\mathbf{x}}(0)=q_{\tilde{\mathbf{x}}}(0)$. Then it is clear from
(\ref{FormOfqxt}) that
\begin{equation}\label{diffqxAndqtx}
\begin{aligned}
 q_{\mathbf{x}}(t)-q_{\tilde{\mathbf{x}}}(t)=&\,q_{\mathbf{x}}(0)\biggl[
e^{\int_0^t P_{\mathbf{x}}(u)du}-e^{\int_0^t
P_{\tilde{\mathbf{x}}}(u)du}\biggl]+\biggl[e^{\int_0^t
P_{\mathbf{x}}(u)du}\int_0^t Q_{\mathbf{x}}(u)e^{-\int_0^u
P_{\mathbf{x}}(v)dv}du\\
&\,-e^{\int_0^t P_{\tilde{\mathbf{x}}}(u)du}\int_0^t
Q_{\tilde{\mathbf{x}}}(u)e^{-\int_0^u
P_{\tilde{\mathbf{x}}}(v)dv}du\biggl].
 \end{aligned}
\end{equation}
 Let $M_1=\sup_{0\le t\le
T}e^{\int_0^t P_{\mathbf{x}}(u)du}\bigvee e^{\int_0^t
P_{\tilde{\mathbf{x}}}(u)du}<\infty$, then for $0\le t\le T$,
\begin{equation}\label{diffexpPxtPtxt}
 |e^{\int_0^t P_{\mathbf{x}}(u)du}-e^{\int_0^t
P_{\tilde{\mathbf{x}}}(u)du}| \le M_1 |\int_0^t
P_{\mathbf{x}}(u)-P_{\tilde{\mathbf{x}}}(u)du |.
\end{equation}
Observe that for
$0\le u\le T$,
\begin{equation}\label{diffPxtPtxtu}
\begin{aligned}
&\,|P_{\mathbf{x}}(u)-P_{\tilde{\mathbf{x}}}(u)|\\
=&\,\biggl|\frac{x_1(u)
\alpha_1(u)\beta_1(u)}{h_{\mathbf{x}}(u)}-\frac{\tilde{x}_1(u)
\alpha_1(u)\beta_1(u)}{h_{\tilde{\mathbf{x}}}(u)} \biggl| \\
=&\,\biggl|
\frac{\alpha_1(u)\alpha_2(u)\beta_1(u)(1+\beta_2(u))[x_1(u)\tilde{x}_2(u)-
\tilde{x}_1(u)x_2(u)]}{h_{\mathbf{x}}(u)h_{\tilde{\mathbf{x}}}(u)}
\\
&\,+\frac{\alpha_1(u)\alpha_3(u)\beta_1(u)(1+\beta_3(u))[x_1(u)\tilde{x}_3(u)-
\tilde{x}_1(u)x_3(u)]}{h_{\mathbf{x}}(u)h_{\tilde{\mathbf{x}}}(u)}\biggl|
\\
=&\,
\biggl|\frac{\alpha_1(u)\alpha_2(u)\beta_1(u)(1+\beta_2(u))[x_1(u)(\tilde{x}_2(u)-x_2(u))+(
x_1(u)-\tilde{x}_1(u))x_2(u)]}{h_{\mathbf{x}}(u)h_{\tilde{\mathbf{x}}}(u)} \\
&\,+\frac{\alpha_1(u)\alpha_3(u)\beta_1(u)(1+\beta_3(u))[x_1(u)(\tilde{x}_3(u)-x_3(u))+(x_1(u)-
\tilde{x}_1(u))x_3(u)]}{h_{\mathbf{x}}(u)h_{\tilde{\mathbf{x}}}(u)}\biggl|
\end{aligned}
\end{equation}
\begin{eqnarray*}
 &\le M_2 \|\mathbf{x}(u)-\tilde{\mathbf{x}}(u)\|,\hspace{10.5cm}
 \end{eqnarray*}
 where $M_2$ depends on $B_L$,
$B_U$ in (\ref{bdhxt}) and $\alpha_i$ and $\beta_i$ ($1\le i\le 3$).
Therefore, for $0\le s\le T$,
\begin{equation}\label{IntegDiffExpPxtPtxt}
\begin{aligned}
\int_0^s |e^{\int_0^t P_{\mathbf{x}}(u)du}-e^{\int_0^t
P_{\tilde{\mathbf{x}}}(u)du}|^2 dt \le &\,\int_0^s M_1^2
\biggl[\int_0^t
|P_{\mathbf{x}}(u)-P_{\tilde{\mathbf{x}}}(u)|du\biggl]^2 dt  \\
\le &\, M_1^2 M_2^2 \int_0^s t \int_0^t
\|\mathbf{x}(u)-\tilde{\mathbf{x}}(u)\|^2 du dt \\
 \le &\, \frac{ M_1^2 M_2^2 T^2}{2} \int_0^s
\|\mathbf{x}(u)-\tilde{\mathbf{x}}(u)\|^2 du.
\end{aligned}
\end{equation}

Next we consider $\int_0^s [e^{\int_0^t P_{\mathbf{x}}(u)du}\int_0^t
Q_{\mathbf{x}}(u)e^{-\int_0^u P_{\mathbf{x}}(v)dv}du-e^{\int_0^t
P_{\tilde{\mathbf{x}}}(u)du}\int_0^t
 Q_{\tilde{\mathbf{x}}}(u)\newline \times e^{-\int_0^u
P_{\tilde{\mathbf{x}}}(v)dv}du]^2 dt $. Notice that
\begin{equation}\label{secondtermDiff}
\begin{aligned}
&\,\biggl|e^{\int_0^t P_{\mathbf{x}}(u)du}\int_0^t
Q_{\mathbf{x}}(u)e^{-\int_0^u P_{\mathbf{x}}(v)dv}du-e^{\int_0^t
P_{\tilde{\mathbf{x}}}(u)du}\int_0^t
Q_{\tilde{\mathbf{x}}}(u)e^{-\int_0^u
P_{\tilde{\mathbf{x}}}(v)dv}du \biggl|  \\
\le&\, \biggl|e^{\int_0^t P_{\mathbf{x}}(u)du}\int_0^t
Q_{\mathbf{x}}(u)e^{-\int_0^u P_{\mathbf{x}}(v)dv}du-e^{\int_0^t
P_{\tilde{\mathbf{x}}}(u)du}\int_0^t Q_{\mathbf{x}}(u)e^{-\int_0^u
P_{\mathbf{x}}(v)dv}du \biggl|\\
&\,+\biggl|e^{\int_0^t P_{\tilde{\mathbf{x}}}(u)du}\int_0^t
Q_{\mathbf{x}}(u)e^{-\int_0^u P_{\mathbf{x}}(v)dv}du-e^{\int_0^t
P_{\tilde{\mathbf{x}}}(u)du}\int_0^t
Q_{\tilde{\mathbf{x}}}(u)e^{-\int_0^u
P_{\mathbf{x}}(v)dv}du\biggl| \\
&\,+\biggl|e^{\int_0^t P_{\tilde{\mathbf{x}}}(u)du}\int_0^t
Q_{\tilde{\mathbf{x}}}(u)e^{-\int_0^u
P_{\mathbf{x}}(v)dv}du-e^{\int_0^t
P_{\tilde{\mathbf{x}}}(u)du}\int_0^t
Q_{\tilde{\mathbf{x}}}(u)e^{-\int_0^u P_{\tilde{\mathbf{x}}}(v)dv}du
\biggl| \\
=&\, I_1(t)+I_2(t)+I_3(t),
 \end{aligned}
 \end{equation}

where $I_i(t)$ denotes the $i$-term on the right hand side of the
inequality in (\ref{secondtermDiff}). Since $P_{\mathbf{x}}(t)$ and
$Q_{\mathbf{x}}(t)$ are continuous functions, there exists $M_3>0$,
such that for any $0\le t\le T$, $|\int_0^t
Q_{\mathbf{x}}(u)e^{-\int_0^u P_{\mathbf{x}}(v)dv}du|\le M_3$. Then
it follows by (\ref{IntegDiffExpPxtPtxt}) that for $0\le s\le T$,
\begin{equation}\label{IntegI1t}
\begin{aligned}
\int_0^s I_1(t)^2dt \le &\, M_3^2 \int_0^s|e^{\int_0^t
P_{\mathbf{x}}(u)du}-e^{\int_0^t
P_{\tilde{\mathbf{x}}}(u)du}|^2 dt \\
 \le &\, \frac{ M_1^2 M_2^2 M_3^2 T^2}{2} \int_0^s
\|\mathbf{x}(u)-\tilde{\mathbf{x}}(u)\|^2 du.
\end{aligned}
\end{equation}

As to $I_2(t)$, we have
\begin{equation}\label{bdI2tExamp}
\begin{aligned}
I_2(t)^2=&\,e^{2\int_0^t
P_{\tilde{\mathbf{x}}}(u)du}\biggl|\int_0^t [Q_{\mathbf{x}}(u)-
Q_{\tilde{\mathbf{x}}}(u)]e^{-\int_0^u
P_{\mathbf{x}}(v)dv}du\biggl|^2 \\
\le&\, M_4^4 t \int_0^t [Q_{\mathbf{x}}(u)-
Q_{\tilde{\mathbf{x}}}(u) ]^2du,
\end{aligned}
\end{equation}
where $M_4>0$ satisfies that for any $0\le t\le T$, $e^{\int_0^t
P_{\tilde{\mathbf{x}}}(u)du}\le M_4$ and $e^{-\int_0^t
P_{\mathbf{x}}(u)du}\le M_4$. Notice that for $0\le u\le T$,
\begin{equation}\label{diffQxtQtxtu}
\begin{aligned}
&\,|Q_{\mathbf{x}}(u)-Q_{\tilde{\mathbf{x}}}(u)|
\\
\le
&\,\biggl|\frac{[\alpha_1(u)\delta_2(u)-\alpha_2(u)\delta_1(u)(1+\beta_2(u))]+\alpha_1(u)\alpha_2(u)\beta_1(u)(1+\beta_2(u))\log
F(u) }{h_{\mathbf{x}}(u)h_{\tilde{\mathbf{x}}}(u)}\biggl| \\
&\, \times |[x_2(u)\tilde{x}_1(u)- x_1(u)\tilde{x}_2(u)]| \\
&\, +\biggl|\frac{[x_3(u)\tilde{x}_1(u)-
x_1(u)\tilde{x}_3(u)][\alpha_1(u)\delta_3(u)-\alpha_3(u)\delta_1(u)(1+\beta_3(u))]}{h_{\mathbf{x}}(u)h_{\tilde{\mathbf{x}}}(u)}\biggl|\\
&\, +\biggl|\frac{[x_3(u)\tilde{x}_2(u)-
x_2(u)\tilde{x}_3(u)][\alpha_2(u)\delta_3(u)(1+\beta_2(u))-\alpha_3(u)\delta_2(u)(1+\beta_3(u))]}{h_{\mathbf{x}}(u)h_{\tilde{\mathbf{x}}}(u)}\biggl|\\
\le &M_5 \|\mathbf{x}(u)-\tilde{\mathbf{x}}(u)\|,
\end{aligned}
\end{equation}
where $M_5$ depends on $B_L$, $B_U$ in (\ref{bdhxt}) and $\alpha_i$,
$\beta_i$, $\delta_i$ ($1\le i\le 3$) and $\log F$. Then it follows
by (\ref{bdI2tExamp}) and (\ref{diffQxtQtxtu}) that for $0\le s\le
T$,
\begin{equation}\label{IntegI2t}
\begin{aligned}
\int_0^s I_2(t)^2dt \le &\, M_4^2 M_5^2 \int_0^s t
\int_0^t \|\mathbf{x}(u)-\tilde{\mathbf{x}}(u)\|^2 du dt \\
 \le &\, \frac{ M_4^4 M_5^2 T^2}{2} \int_0^s
\|\mathbf{x}(u)-\tilde{\mathbf{x}}(u)\|^2 du.
\end{aligned}
\end{equation}

As to $I_3(t)$, we have
\begin{equation}\label{bdI3tExamp}
\begin{aligned}
I_3(t)^2=&\,e^{2\int_0^t P_{\tilde{\mathbf{x}}}(u)du}\biggl|\int_0^t
Q_{\tilde{\mathbf{x}}}(u)[e^{-\int_0^u
P_{\mathbf{x}}(v)dv}-e^{-\int_0^u
P_{\tilde{\mathbf{x}}}(v)dv}]du\biggl|^2 \\
\le&\, M_4^2 M_6^2 t \int_0^t [e^{-\int_0^u
P_{\mathbf{x}}(v)dv}-e^{-\int_0^u P_{\tilde{\mathbf{x}}}(v)dv}]^2
du,
\end{aligned}
\end{equation}
where $M_6>0$ satisfies that for any $0\le t\le T$,
$|Q_{\tilde{\mathbf{x}}}(t)|\le M_6$. Similar to
(\ref{IntegDiffExpPxtPtxt}), we can prove that there exists $M_7>0$,
for any $0\le t\le T$, such that
\begin{equation}\label{I3tM7}
\int_0^t [e^{-\int_0^u P_{\mathbf{x}}(v)dv}-e^{-\int_0^u
P_{\tilde{\mathbf{x}}}(v)dv}]^2 du\le M_7^2  \int_0^t
\|\mathbf{x}(u)-\tilde{\mathbf{x}}(u)\|^2 du.
\end{equation}
 Then it follows that
for $0\le s\le T$,
\begin{equation}\label{IntegI3t}
\int_0^s I_3(t)^2dt \le \frac{ M_4^2
M_6^2 M_7^2 T^2}{2} \int_0^s
\|\mathbf{x}(u)-\tilde{\mathbf{x}}(u)\|^2 du.
\end{equation}

Thus, by (\ref{secondtermDiff}), (\ref{IntegI1t}), (\ref{IntegI2t}),
(\ref{IntegI3t}), we get
\begin{equation}\label{SecdondDiffInt}
\begin{aligned}
&\,\int_0^s \biggl|e^{\int_0^t P_{\mathbf{x}}(u)du}\int_0^t
Q_{\mathbf{x}}(u)e^{-\int_0^u P_{\mathbf{x}}(v)dv}du-e^{\int_0^t
P_{\tilde{\mathbf{x}}}(u)du}\int_0^t
Q_{\tilde{\mathbf{x}}}(u)e^{-\int_0^u P_{\tilde{\mathbf{x}}}(v)dv}du
\biggl|^2 dt \\
\le &\,\frac{3(M_1^2 M_2^2 M_3^2 T^2+ M_4^4 M_5^2 T^2+M_4^2 M_6^2
M_7^2 T^2)}{2} \int_0^s \|\mathbf{x}(u)-\tilde{\mathbf{x}}(u)\|^2
du.
\end{aligned}
\end{equation}

By (\ref{diffqxAndqtx}), (\ref{IntegDiffExpPxtPtxt}) and
(\ref{SecdondDiffInt}), (\ref{DiffOfTwoSolnsforqExamp}) is true for
some $M>0$.

The uniqueness of the $C_{K^3\times R}[0,\infty)$-martingale problem
for $(G_A, \mu)$ follows by (\ref{DiffOfTwoSolnsforqExamp}) and
Lemma \ref{lem:UniqForIAFFM}. \qed

\section*{Appendix B: Proof of Lemma \ref{rem:verifyCondiEofCor3.6Chp2}}

\proof Let $f\in C^2_c(K\times R)$ and fix $T>0$. Assume that
$\tilde{K}$ is the support of $f$. For $0\le t\le T$, and ${\mathbf
y}=(\mathbf{x}, q)\in K_N\times F_N$, by Taylor's expansion,
(\ref{eq:meanIAFFM}) and (\ref{eq:RecLogPriceqNtgN})
\begin{equation}\label{eq:MarkTayGIAFFM}
\begin{aligned}&\,[S_{N,[Nt]}-I]f({\mathbf y})\\
 =&\,E[(\mathbf{X}^N(t+\frac{1}{N})-\mathbf{x})\frac{\partial f}{\partial {\mathbf x}}({\mathbf y})'|\mathbf{Y}^N(t)={\mathbf
 y}]+E[(q^N(t+\frac{1}{N})-q)\frac{\partial f}{\partial q}({\mathbf
 y})|\mathbf{Y}^N(t)={\mathbf y}] \\
  &\,+E[(\mathbf{Y}^N(t+\frac{1}{N})-{\mathbf y})\frac{\partial^2 f}{\partial {\mathbf y}^2}(\mathbf{Y}^N_*(t))(\mathbf{Y}^N(t+\frac{1}{N})-{\mathbf y})'|\mathbf{Y}^N(t)={\mathbf y}] \\
 =&\,\frac{1}{N}\mathbf{x}A_N(\frac{[Nt]}{N},{\mathbf y})\frac{\partial f}{\partial {\mathbf x}}({\mathbf
 y})'+\frac{1}{N}E\biggl[g_N\biggl(\frac{[Nt]+1}{N},\mathbf{X}^N(t+\frac{1}{N}),q\biggl)\frac{\partial f}{\partial q}({\mathbf
 y})\biggl|\mathbf{Y}^N(t)={\mathbf y}\biggl] \\
 &\,+ E[(\mathbf{Y}^N(t+\frac{1}{N})-{\mathbf y})\frac{\partial^2 f}{\partial {\mathbf y}^2}(\mathbf{Y}^N_*(t))(\mathbf{Y}^N(t+\frac{1}{N})-{\mathbf y})'|\mathbf{Y}^N(t)={\mathbf y}],
\end{aligned}
\end{equation}
 where $\mathbf{Y}^N_*(t)={\mathbf y}+\theta^N_t
(\mathbf{Y}^N(t+\frac{1}{N})-{\mathbf y})$, for some $\theta^N_t\in
(0,1)$.

Notice that
\begin{equation}\label{gNtgtreformIAFFM}
\begin{aligned}
&\,g_N\biggl(\frac{[Nt]+1}{N},\mathbf{X}^N(t+\frac{1}{N}),q\biggl)-g(t,\mathbf{x},q)
\\ =&\,\biggl[
g_N\biggl(\frac{[Nt]+1}{N},\mathbf{X}^N(t+\frac{1}{N}),q\biggl)-g\biggl(\frac{[Nt]+1}{N},\mathbf{X}^N(t+\frac{1}{N}),q\biggl)\biggl]\\
&\,+\biggl[g\biggl(\frac{[Nt]+1}{N},\mathbf{X}^N(t+\frac{1}{N}),q\biggl)-g\biggl(\frac{[Nt]+1}{N},\mathbf{x},q\biggl)\biggl]\\
&\,+\biggl[g\biggl(\frac{[Nt]+1}{N},\mathbf{x},q\biggl)-g(t,\mathbf{x},q)\biggl]
\end{aligned}
\end{equation}
and
\begin{equation}\label{gNtgtIAFFM}
\begin{aligned}&\,g\biggl(\frac{[Nt]+1}{N},\mathbf{X}^N(t+\frac{1}{N}),q\biggl)-g\biggl(\frac{[Nt]+1}{N},\mathbf{x},q\biggl)
\\
=&\, (\mathbf{X}^N(t+\frac{1}{N})-\mathbf{x})\frac{\partial
g}{\partial {\mathbf
x}}(\frac{[Nt]+1}{N},\mathbf{X}^N_{\star}(t+\frac{1}{N}),q)',
\end{aligned}
\end{equation} where
$\frac{\partial \mathbf{g}}{\partial \mathbf{x}}=[\frac{\partial
g}{\partial x_1}, \cdots, \frac{\partial g}{\partial x_r}]$,
$\mathbf{X}^N_{\star}(t+\frac{1}{N})=\mathbf{x}+\zeta^N_{t+\frac{1}{N}}
(\mathbf{X}^N(t+\frac{1}{N})-\mathbf{x})$, for some
$\zeta^N_{t+\frac{1}{N}}\in (0,1)$.

Then it follows that
\begin{equation}\label{eq:NSNNtGAtIAFFM}
\begin{aligned}{}&|N[S_{N,[Nt]}-I]f({\mathbf y})-G_A(t)f({\mathbf y})| \\
 \le&\,\biggl|\mathbf{x}[A_N(\frac{[Nt]}{N},{\mathbf y})-A(t,{\mathbf y})]\frac{\partial f}{\partial {\mathbf x}}({\mathbf
 y})'\biggl|+E\biggl[\biggl|\biggl(g_N\biggl(\frac{[Nt]+1}{N},\mathbf{X}^N(t+\frac{1}{N}),q\biggl)\\
 &\,-g\biggl(\frac{[Nt]+1}{N},\mathbf{X}^N(t+\frac{1}{N}),q\biggl)\biggl)
\frac{\partial f}{\partial q}({\mathbf
y})\biggl|\biggl|\mathbf{Y}^N(t)={\mathbf y}\biggl]
\\
&\,+E\biggl[\biggl|(\mathbf{X}^N(t+\frac{1}{N})-\mathbf{x})\frac{\partial
g}{\partial {\mathbf
x}}(\frac{[Nt]+1}{N},\mathbf{X}^N_{\star}(t+\frac{1}{N}),q)'\frac{\partial
f}{\partial
q}({\mathbf y})\biggl|\biggl|\mathbf{Y}^N(t)={\mathbf y}\biggl] \\
&\,+E\biggl[\biggl|\biggl(g\biggl(\frac{[Nt]+1}{N},\mathbf{x},q\biggl)-g(t,\mathbf{x},q)\biggl)\frac{\partial
f}{\partial q}({\mathbf
y})\biggl|\biggl|\mathbf{Y}^N(t)={\mathbf y}\biggl] \\
 &\,+ NE[\biggl|(\mathbf{Y}^N(t+\frac{1}{N})-{\mathbf y})\frac{\partial^2
f}{\partial {\mathbf
y}^2}(\mathbf{Y}^N_*(t))(\mathbf{Y}^N(t+\frac{1}{N})-{\mathbf
y})'\biggl|\biggl|\mathbf{Y}^N(t)={\mathbf y}] \\
=&\, I_1(N,{\mathbf y},t)+I_2(N,{\mathbf y},t)+I_3(N,{\mathbf
y},t)+I_4(N,{\mathbf y},t)+I_5(N,{\mathbf y},t).
\end{aligned}
\end{equation}

 Notice that for each
$(\mathbf{x},q)\in K_N\times R$, $A_N(t,\mathbf{x},q)$ is a constant
on $[\frac{k}{N}, \frac{k+1}{N})$ for each $k\ge 0$, we have
$A_N(t,\mathbf{x},q)=A_N(\frac{[Nt]}{N},\mathbf{x},q)$ for $t\ge 0$.
The norm of matrices is defined as follows: for real-valued matrix
$B=(b_{i,j})_{r\times r}$, $\|B\|=\sum_{i,j=1}^r|b_{i,j}|$. The norm
of real-valued $r$-dimensional vector ${\mathbf
y}=(y_1,\cdots,y_r)'$ is the Euclidean norm $\|{\mathbf y}\|$. It
follows that $\|B {\mathbf y}\|\le \|B\| \cdot\|\mathbf{y}\|$. Then
we have
\begin{equation}\label{eq:I1NytIAFFM}
\begin{aligned}I_1(N,{\mathbf y},t) \le
&\,\|A_N(t,\mathbf{y})-A(t,\mathbf{y})\|\cdot \|\frac{\partial
f}{\partial {\mathbf x}}({\mathbf y})'\|.
\end{aligned}
\end{equation}
Then it follows by (\ref{eq:ANConvAUniOnF}) that
\begin{equation}\label{eq:I1NytUniformConIAFFM}
\begin{aligned}\lim_{N\to \infty}\sup_{0\le t\le
T}\sup_{q\in F_N}\sup_{\mathbf{x}\in K_N} I_1(N,\mathbf{x},q,t)=0.
\end{aligned}
\end{equation}

Next, we consider $I_5(N,{\mathbf y},t)$. Let $\|\frac{\partial^2
f}{\partial {\mathbf y}^2}\|=\max_{1\le i,j\le
r+1}\|\frac{\partial^2 f}{\partial y_i\partial y_j}\|$ and define
$g_{N,\tilde{K}}$ as follows: $
g_{N,\tilde{K}}(t,\mathbf{y})=g_N(t,\mathbf{y})$, if
$(t,\mathbf{y})\in [0,\infty)\times (\tilde{K}\bigcap (K_N\times
R))$; $g_{N,\tilde{K}}(t,\mathbf{y})=0$ otherwise. By H\"{o}lder
inequality and Cauchy-Schwartz inequality,
\begin{equation}\label{eq:I5NytIAFFM}
\begin{aligned}
&\,I_5(N,{\mathbf y},t) \\
\le &\,\|\frac{\partial^2 f}{\partial {\mathbf
y}^2}\|N\biggl\{\sum_{i=1}^{r}\sum_{j=1}^{r}E[|(X^N_i(t+\frac{1}{N})-x_i)(X^N_j(t+\frac{1}{N})-x_j)|\biggl|\mathbf{Y}^N(t)={\mathbf
y}] \\
&\,+\frac{2}{N}\sum_{i=1}^r
E\biggl[|(X^N_i(t+\frac{1}{N})-x_i)g_{N,\tilde{K}}\biggl(\frac{[Nt]+1}{N},\mathbf{X}^N(t+\frac{1}{N}),q\biggl)|\biggl|\mathbf{Y}^N(t)={\mathbf
y}\biggl]    \\
&\,+\frac{1}{N^2}E\biggl[|g_{N,\tilde{K}}\biggl(\frac{[Nt]+1}{N},\mathbf{X}^N(t+\frac{1}{N}),q\biggl)|^2\biggl|\mathbf{Y}^N(t)={\mathbf
y}\biggl]\biggl\} \\
\le &\,\|\frac{\partial^2 f}{\partial {\mathbf
y}^2}\|N\biggl\{\sum_{i=1}^{r}\sum_{j=1}^{r}\biggl(E[(X^N_i(t+\frac{1}{N})-x_i)^2\biggl|\mathbf{Y}^N(t)={\mathbf
y}]\\
&\,\times
E[(X^N_j(t+\frac{1}{N})-x_j)^2\biggl|\mathbf{Y}^N(t)={\mathbf
y}]\biggl)^{\frac{1}{2}}  \\
&\,+\frac{2}{N}\|g_{N,\tilde{K}}\|_{T+1}\sum_{i=1}^r
\biggl(E\biggl[(X^N_i(t+\frac{1}{N})-x_i)^2\biggl|\mathbf{Y}^N(t)={\mathbf
y}\biggl] \biggl)^{\frac{1}{2}}
+\frac{1}{N^2}\|g_{N,\tilde{K}}\|_{T+1}^2 \biggl\},
\end{aligned}
\end{equation}
where $\|g_{N,\tilde{K}}\|_{T+1}=\sup_{0\le t\le T+1}\sup_{{\mathbf
y} \in \tilde{K}\bigcap (K_N\times R)}|g_N(t,y)|$. By
(\ref{VarphiNconvToVarphi}), (\ref{PsiNconvToPsi}) and the
assumption that $\varphi(t, \mathbf{x}, q)$ and $\psi(t, \mathbf{x},
q)$ are continuous on $[0,T]\times K\times R$ for any $T>0$, it
follows that
\begin{equation}\label{gNtdKbd}
\sup_{N\ge 1}\|g_{N,\tilde{K}}\|_{T+1}<\infty.
\end{equation}
 For fixed $1\le i\le
r$, by (\ref{eq:squareOfDifIAFFM}), we get
\begin{equation}\label{eq:CalsquareOfDifIAFFM}
\begin{aligned}{}& N^2E[(X^N_i(t+\frac{1}{N})- x_i)^2|\mathbf{Y}^N(t)={\mathbf y}] \\
  =&\, \mathbf{x}A_{N,\cdot,i}(\frac{[Nt]}{N},\mathbf{y})+(\mathbf{x}A_{N,\cdot,i}(\frac{[Nt]}{N},\mathbf{y}))^2
  -2x_i a_{N,i,i}(\frac{[Nt]}{N},\mathbf{y})- \sum_{k=1}^r \frac{x_k}{N} a_{N,k,i}^2 (\frac{[Nt]}{N},\mathbf{y}) \\
  \le &\, \sum_{l=1}^r  |A_{N,l,i}(t,\mathbf{y})|+(\sum_{l=1}^r  |A_{N,l,i}(t,\mathbf{y})|)^2+2 |a_{N,i,i}(t,\mathbf{y})|.
  \end{aligned}
\end{equation}
Then it follows by (\ref{eq:I5NytIAFFM}), (\ref{gNtdKbd}),
(\ref{eq:CalsquareOfDifIAFFM}), (\ref{eq:ANConvAUniOnF}) and the
assumption that $A(t,\mathbf{x},q)$ is bounded on $[0, T]\times
K\times R$, that
\begin{equation}\label{eq:I5NytUniformConIAFFM}
\begin{aligned}
\lim_{N\to\infty}\sup_{0\le t\le T}\sup_{q\in
F_N}\sup_{\mathbf{x}\in K_N} I_5(N,\mathbf{x},q,t)=0.
  \end{aligned}
\end{equation}

Similarly,  by (\ref{VarphiNconvToVarphi}) and
(\ref{PsiNconvToPsi}), we can get
\begin{equation}\label{eq:I2NytUniformConIAFFM}
\begin{aligned}
\lim_{N\to \infty}\sup_{0\le t\le T}\sup_{q\in
F_N}\sup_{\mathbf{x}\in K_N} I_2(N,\mathbf{x},q,t)=0,
  \end{aligned}
\end{equation}
and by (\ref{eq:CalsquareOfDifIAFFM}), (\ref{eq:ANConvAUniOnF}) and
some calculations, we can get that
\begin{equation}\label{eq:I3NytUniformConIAFFM}
\begin{aligned}
\lim_{N\to \infty}\sup_{0\le t\le T}\sup_{q\in
F_N}\sup_{\mathbf{x}\in K_N} I_3(N,\mathbf{x},q,t)=0.
  \end{aligned}
\end{equation}

Notice that for any $(\mathbf{x},q)\in K\times R$ and $0\le t \le
T$,
\begin{equation}\label{I4NtBd}
I_4(N,\mathbf{x},q,t)\le \sup_{0\le t\le T}\sup_{(\mathbf{x},q)\in
\tilde{K}}
\biggl|g\biggl(\frac{[Nt]+1}{N},\mathbf{x},q\biggl)-g(t,\mathbf{x},q)\biggl|\cdot\sup_{\mathbf{y}\in
\tilde{K}}|\frac{\partial f}{\partial q}({\mathbf y})|,
\end{equation}
 then it follows by the uniform continuity of $g(t,\mathbf{x},q)$ on the
compact set $[0,T+1]\times \tilde{K}$ that
\begin{equation}\label{eq:I4NytUniformConIAFFM}
\begin{aligned}
\lim_{N\to \infty}\sup_{0\le t\le T}\sup_{q\in
F_N}\sup_{\mathbf{x}\in K_N} I_4(N,\mathbf{x},q,t)=0,
  \end{aligned}
\end{equation}
and (\ref{GNfGfIAFFM}) is proved. \qed

\section*{Appendix C: Proof of Lemma \ref{lem:LimitPtBeingCont}}

Assume that $\mathbf{Y}=(\mathbf{X}, Q)$, where
$\mathbf{X}=(X_1,\cdots, X_r)$, is a limit point of $\mathbf{Y}^N$.
By Corollary \ref{cor:YNtight}, $\mathbf{Y}$ is a solution of the
$D_{K\times R}[0,\infty)$-martingale problem for $(G_A, \mu)$
restricted on $C^2_c(K\times R)$. We just need to show that
$\mathbf{Y}$ is continuous almost surely.

Fix $T>0$. Since $A(t, \mathbf{x}, q)$, $\varphi(t,\mathbf{x},q)$
and $\psi(t,\mathbf{x},q)$ are bounded on $[0, T]\times K\times R$,
there exists $C_T>0$, such that for $1\le i\le r$,
\begin{equation}\label{bdforBiVarphiPsi}
\sup_{0\le t\le T}\sup_{\mathbf{y}\in K\times
R}|b_i(t,\mathbf{y})|\le C_T, \mbox{ }\sup_{0\le t\le
T}\sup_{\mathbf{y}\in K\times R}|\varphi(t,\mathbf{y})|\le C_T,
\mbox{ }\sup_{0\le t\le T}\sup_{\mathbf{y}\in K\times
R}|\psi(t,\mathbf{y})|\le C_T.
\end{equation}
By (\ref{VarphiNconvToVarphi}) and (\ref{PsiNconvToPsi}), there
exists $N_0$ such that for $N>N_0$, we have
\begin{equation}\label{bdforVarphiNPsiN}
 \mbox{ }\sup_{0\le t\le
T}\sup_{\mathbf{y}\in K_N\times R}|\varphi_N(t,\mathbf{y})|\le
C_T+1, \mbox{ }\sup_{0\le t\le T}\sup_{\mathbf{y}\in K_N\times
R}|\psi_N(t,\mathbf{y})|\le C_T+1.
\end{equation}

Let $f_i(\mathbf{x},q)=x_i$, $1\le i\le r$. It follows that
$f_i^j\in C^2_c(K\times R)$ for $1\le i\le r$ and $1\le j\le 4$. We
assume that $\{\mathscr{F}_s, 0\le s<\infty\}$ is the filtration to
which the $D_{K\times R}[0,\infty)$-martingale problem for $(G_A,
\mu)$ referred. Then $X^j_i(t)-j\int_0^t (X_i(u))^{j-1}
\mathbf{X}(u) \times A_{\cdot, i}(u, \mathbf{Y}(u))du
=X^j_i(t)-j\int_0^t (X_i(u))^{j-1} b_i(u, \mathbf{Y}(u))du$ is an
$\{\mathscr{F}_s\}$-martingale for $1\le i\le r$ and $1\le j\le 4$.
Let $0\le s<t\le T$, and fix $1\le i\le r$, it follows that
\begin{equation}\label{eq:DiffToPowerj}
E[(X_i(t)-X_i(s))^j]=jE[\int_s^t(X_i(u)-X_i(s))^{j-1}b_i(u,
\mathbf{Y}(u))du]
\end{equation}
 for $1\le j\le 4$. Then by
(\ref{eq:DiffToPowerj}), for $j=2$, we have
\begin{equation}\label{eq:DiffToPower2}
\begin{aligned}
 E[(X_i(t)-X_i(s))^2]\le &\, 2 E[\int_s^t |X_i(u)-X_i(s)|\,|b_i(u, \mathbf{Y}(u))|du] \\
\le & \, 2 C_T(t-s).
\end{aligned}
\end{equation}
 By (\ref{eq:DiffToPowerj}) and (\ref{eq:DiffToPower2}), for $j=4$,
\begin{equation}\label{eq:DiffToPower4}
\begin{aligned}
E[(X_i(t)-X_i(s))^4]\le &\, 4 C_T E[\int_s^t |X_i(u)-X_i(s)|^ 3 du] \\
\le &\, 4 C_T E[\int_s^t |X_i(u)-X_i(s)|^ 2 du] \\
\le &\, 4 C_T \int_s^t 2 C_T(u-s) du \\
\le &\, 4 (C_T)^2 (t-s)^2.
\end{aligned}
\end{equation}
 Then by Kolmogorov's Criterion, we
proved for $1\le i\le r$ that $X_i$ is continuous almost surely.

Next, we prove that $Q$ is also continuous almost surely. We
introduce the notations in Chapter 3, Section 10 \cite{Kurtz}. Let
$(E, r)$ be a metric space. For $x\in D_{E}[0, \infty)$, define
\begin{equation}\label{defnJ(x)inEK}
 J(x)=\int_0^\infty e^{-u}[J(x,u)\wedge 1]du,
 \end{equation}
where $J(x,u)$ is defined by
\begin{equation}\label{defnJ(x,u)inEK}
J(x,u)=\sup_{0\le t\le u} r(x(t),x(t-)).
\end{equation}
 Since $\mathbf{Y}$ is a
limit point of $\mathbf{Y}^N$, we have that a subsequence $\{q^{N_k}
\}$ of $\{q^N \}$ converges weakly to $Q$. To prove that $Q$ is
continuous almost surely, by Chapter 3, Theorem 10.2 \cite{Kurtz},
it suffices to prove that $J(q^{N_k})\Rightarrow 0$ as $k\rightarrow
\infty$. In this case $E=R$ and $r$ is the Euclidean metric. It is
enough to show that $\lim_{k\to \infty}E[J(q^{N_k})]=0$.

By Lemma \ref{lem:CompactContainForYN}, $\{q^{N_k} \}$ satisfies the
compact containment condition, i.e. for any $\eta>0$ and $T>0$,
there exists $B_T>0$, such that
\begin{equation}\label{eq:compatCont}
 \inf_{k}P\{|q^{N_k}(t)|\le
B_T, 0\le t\le T\}\ge 1-\eta.
\end{equation}
By the construction of $q^N$, (\ref{eq:RecLogPriceTdqNgN}) and
(\ref{bdforVarphiNPsiN}), for fixed $N_k\ge N_0$ and $0\le u\le T$,
\begin{equation}\label{J(qNk,u)IFMM}
\begin{aligned}
J(q^{N_k},u)=&\,\max_{0\le j\le \frac{[N_k
u]}{N_k}}|\tilde{q}^{N_k}(\frac{j}{N_k})-\tilde{q}^{N_k}(\frac{j-1}{N_k})|
\\
\le &\, \frac{1}{N_k} \max_{0\le j\le \frac{[N_k u]}{N_k}} \biggl[|
\varphi_{N_k}(\frac{j}{N_k},\mathbf{X}(\frac{j}{N_k}),\tilde{q}^{N_k}(\frac{j-1}{N_k}))\tilde{q}^{N_k}(\frac{j-1}{N_k})|
\\
&\, +|\psi_{N_k}(\frac{j}{N_k},\mathbf{X}(\frac{j}{N_k}),
\tilde{q}^{N_k}(\frac{j-1}{N_k}))|\biggl]
\\
\le &\, \frac{1}{N_k} (C_T+1)(B_T+1),
\end{aligned}
\end{equation}
 on the event $F_{k,T}=\{|q^{N_k}(t)|\le B_T, 0\le t\le T\}$. It follows by
 (\ref{eq:compatCont}) and (\ref{J(qNk,u)IFMM}) that
\begin{equation}\label{eq:ExpectedValueJqNk}
\begin{aligned}
E[J(q^{N_k})]\le&\, e^{-T}+E[\int_0^T
 e^{-u}(J(q^{N_k},u)\wedge 1)du ] \\
 \le&\, e^{-T}+E[\int_0^T
 e^{-u}(J(q^{N_k},u)\wedge 1)du \chi_{F_{k,T}}]+P(F_{k,T}^c) \\
 \le&\, e^{-T}+ \frac{T}{N_k}
 (C_T+1)(B_T+1)+\eta.
   \end{aligned}
   \end{equation}
Let $k\rightarrow \infty$ and then let $T\rightarrow \infty$,
$\eta\rightarrow 0$, we proved that $\lim_{k\to
\infty}E[J(q^{N_k})]=0$. \qed

\section*{Acknowledgements}

The author would like to thank Professor Donald A. Dawson for his
valuable suggestions of the references, helpful discussions, and the
enlightening comments. The author also thanks Dr. Ulrich Horst for
sending his ph.D Thesis and other papers, and especially for his
comment on the economic interpretation of the zeros of the
denominator in the function $g$.

\end{document}